\documentclass[10pt]{article}
 \usepackage[top=1in,bottom=1in,left=1.5in,right=1.5in]{geometry}
  \usepackage{amsmath,amssymb}
  \usepackage{latexsym}
 \usepackage[dvips]{pstcol} 
  \usepackage{pst-node}
  \usepackage[dvips]{graphicx}

  \newcommand{\C}{\mathbb{C}}
  \newcommand{\N}{\mathbb{N}}
  \newcommand{\R}{\mathbb{R}}
  \newcommand{\Z}{\mathbb{Z}}

  \newcommand{\bd}{\mathbf{d}}
  \newcommand{\e}{\mathbf{e}}

  \newcommand{\p}{\mathbf{p}}
  \newcommand{\bP}{\mathbb{P}}
  
  \newcommand{\q}{\mathbf{q}}

  \renewcommand{\u}{\mathbf{u}}
  \renewcommand{\v}{\mathbf{v}}
  
  \newcommand{\w}{\mathbf{w}}
  \newcommand{\x}{\mathbf{x}}
  \newcommand{\y}{\mathbf{y}}
  \newcommand{\z}{\mathbf{z}}
  \newcommand{\0}{\mathbf{0}}
  \newcommand{\1}{\mathbf{1}}

  \newcommand{\Cp}{\mathrm{Cap\;}}
  \newcommand{\lan}{\langle}
  \newcommand{\ran}{\rangle}
  \newcommand{\an}[1]{\lan#1\ran}
  \def\diag{\mathop{{\rm diag}}\nolimits}
  \newcommand{\hs}{\hspace*{\parindent}}
  \newcommand{\proof}{\hs \textbf{Proof.\ }}

  \newcommand{\trans}{^\top}
  \newcommand{\qed}{\hspace*{\fill} $\Box$\\}
  \newcommand{\per}{\mathop{\mathrm{perm}}\nolimits}
  \newcommand{\haff}{\mathrm{haf\;}}

  \newcommand{\rS}{\mathrm{S}}

  \newcommand{\perm}{\mathrm{perm\;}}
  
  \newcommand{\set}[1]{\{#1\}}

  \newtheorem{theo}{\bfseries \hs Theorem}[section]
  \newtheorem{defn}[theo]{\bfseries \hs Definition}
  \newtheorem{prop}[theo]{\bfseries \hs Proposition}
  \newtheorem{lemma}[theo]{\bfseries \hs Lemma}
  \newtheorem{corol}[theo]{\bfseries \hs Corollary}
  \newtheorem{con}[theo]{\bfseries \hs Conjecture}
  
  \newtheorem{rem}[theo]{\bfseries \hs Remark}
  \numberwithin{equation}{section} 

\begin{document}

 \title{Generalized Friedland-Tverberg inequality:\\ applications
 and extensions}

 \author{Shmuel Friedland
 \\Department of Mathematics, \\ Statistics and Computer
 Science\\ University of Illinois at Chicago
 \\ Chicago, Illinois 60607-7045\\
 friedlan@uic.edu
 \and Leonid Gurvits
 \\Los Alamos\\National Laboratories \\
 Los Alamos\\NM 87545\\
 gurvits@lanl.gov}

 \date{August 21, 2006}

 \maketitle

 \begin{abstract}
 We derive here the Friedland-Tverberg inequality for positive hyperbolic
 polynomials.  This inequality is applied to give lower bounds for
 the number of matchings in $r$-regular bipartite graphs.  It is
 shown that some of these bounds are asymptotically sharp.  We
 improve the known lower bound for the three dimensional monomer-dimer
 entropy. We present Ryser-like formulas for computations
 of matchings in bipartite and general graphs.
 Additional algorithmic applications are given.
     \\[\baselineskip] 2000 Mathematics Subject
     Classification: 05A15, 05A16, 05C70, 05C80, 82B20

 \end{abstract}

 \section{Introduction}

 The aim of this paper is to explore the connections
 between the problem of counting the number of partial matchings
 in graphs and positive hyperbolic polynomials, and Ryser-like
 formulas for partial matchings.
 Given a graph $G=(V,E)$ on $N$ vertices, i.e. $\#V=N$,
 we want to compute the number of $m$-matching,
 i.e. the number of subsets $M$ of edges $E$, where $\#M=m$,
 and no two edges in $M$ have a common vertex.

 Our main results are for bipartite graphs
 $G:=(V_1\cup V_2,E)$, where $E\subset V_1\times V_2$
 and $n=\#V_1=\#V_2$.  Let $A(G)\in \{0,1\}^{n\times n}$
 be the incidence matrix of the bipartite graph $G$.
 Then the number of $m$-matchings in $G$ is equal to $\per_m A(G)$,
 where $\per_m A$ is the sum of $m\times m$ minors of $A\in
 \R^{n\times n}$.  For $m=n$, $\per A(G)$, the permanent
 of $A(G)$ is the number of perfect matchings in $G$.

 It is well known that the computation of the number of perfect
 matching in a general bipartite graph is
 $\#P-complete$.  See \cite{Val} for the first
 proof and \cite{BS} for a simplified proof.

 Ryser's algorithm to compute the permanent of any
  $A\in \R^{n\times n}$ \cite{Rys} remains
 the most efficient exact algorithm, even though it uses
 around $n2^n$ operations.  One can speed up significantly the
 approximate computation of $\per A$.
 One knows to compute the permanent of a nonnegative
 matrix $A=[a_{ij}]\in \R_+^{n\times n}$
 within a simply exponential factor \cite{Bar}.  In the case
 that all the entries of the matrix are uniformly bounded below and above
 by positive constants, one can improve the estimates of the
 exponential errors \cite{LSW} and \cite{FRZ}.
 A fully randomized polynomial approximation scheme \emph{frpas}
 for the number of perfect matchings in a bipartite $G$, and more
 generally for the permanent of a nonnegative matrix is given
 in \cite{JSV}.  This result is generalized in
 \cite{FL} to the number of $m$-matchings in bipartite $G$, and more generally
 to $\per_m A$, for any $A\in \R_+^{l\times n}$.

 We now describe our main results for $\per_m A$, where $A$ is doubly
 stochastic,
 and their applications to lower bounds on partial matchings in
 bipartite graphs.
  Recall that the minimum of the permanent of $n\times n$
 doubly stochastic matrices, denoted by $\Omega_n$, is achieved
 only for the flat matrix $J_n$, whose all entries equal to
 $\frac{1}{n}$.  Thus $\perm B\ge \perm J_n=\frac{n!}{n^n}$ for any
 $B\in \Omega_n$ and this inequality was conjectured by van der
 Waerden \cite{vdW}.  This conjecture was independently proved by
 Egorichev \cite{Ego} and Falikman \cite{Fal}.  We call the above
 inequality Egorichev-Falikman-van der Waerden  (EFW) inequality.
 The asymptotic
 behavior of EFW inequaity is captured by the
 inequality $\perm B\ge e^{-n}$ for any $B\in\Omega_n$.
 This inequality was shown by the first name author \cite{Fri}
 three years before \cite{Ego, Fal}.  Let $\Gamma(n,r)$ be the set
 of all $r$-regular bipartite graphs $G$ on $2n$ vertices.
 For $G\in \Gamma(n,r)$ the matrix $B:=\frac{1}{r}A(G)$
 is doubly stochastic.  Hence the number of perfect matchings
 in $G$ is at least $(\frac{r}{e})^n$.  Thus for $r\ge 3$,
 the number of perfect matchings in $r$-regular bipartite
 graphs grows exponentially, which proves a conjecture by
 Erdos-Renyi \cite{ER}.
 Schrijver \cite{Sch} improved EFW inequality to $r$-regular
 bipartite graphs, whose asymptotic growth is best possible.
 Recently, the second name author
 \cite{Gur1} improved Schrijver's inequality.  Moreover,
 the proof in \cite{Gur1} is significantly simpler and
 transparent.  One of the main tools in the proof in \cite{Gur1}
 is the use of the classical theory of hyperbolic polynomials.

 It was shown by the first name author that $\perm_m
 A\ge \perm_m J_n$ for any $A\in\Omega_n$, and for $m\in [2,n]$
 equality holds only if and only if $A=J_n$ \cite{Fr2}.
 ($\perm_1 A=n$ for each $A\in\Omega_n$.)  This was fact was
 conjectured by Tverberg \cite{Tve}, and is called in this paper the
 Friedland-Tverberg (FT) inequality.  FT inequality gives a lower
 bound on the number of partial matchings in any $G\in\Gamma(n,r)$.

 We derive here the  Schrijver type inequalities
 for $m$ matchings in $r$-regular bipartite graphs on $2n$ vertices.
 This is done using the results and techniques of \cite{Gur1}.
 In particular we give a generalized versions of FT inequality
 to positive homogeneous hyperbolic polynomials, which are of
 independent interest.

 The notion of partial matching in $\Gamma(n,r)$
 can be extended to asymptotic matchings as $n\to\infty$
 as follows.  Given a sequence of $G_l\in\Gamma(n_l,r)$ we can
 consider the quantities
 \begin{eqnarray}\label{defpasmat}
 f(p,\{G_l\}):=\liminf_{n_l\to\infty} \frac{\log \per _{m_l} A(G_l)}{2n_l}, \quad
 F(p,\{G_l\}):=\limsup_{n_l\to\infty} \frac{\log \per _{m_l} A(G_l)}{2n_l},\\
 \textrm{where } n_l\to\infty \textrm{ and } \lim_{l\to\infty} \frac{m_l}{n_l}
 =p\in [0,1].\nonumber
 \end{eqnarray}
 $f(p,\{G_l\})$ and $F(p,\{G_l\})$ can be viewed as the minimal
 and the maximal exponential growth of matchings of density $p$
 of the sequence $G_l,l\in\N$.

 Consider the following special case of the above example.
 Let $C_m$ be a cycle on $m$ vertices.  Note that $C_m$ is bipartite
 if and only if $m$ s even.  Fix a positive integer $d$ and let
 $T_{2l,d}:=\underbrace{C_{2l}\times \ldots\times C_{2l}}_d$, be
 bipartite toroidal grid on $(2l)^d$ vertices. Note that $T_{2l,d}\in
 \Gamma(2^{d-1}l^d, 2d)$.  It is shown in \cite{Ha1} that
 $f(p,\{T_{2l,d}\})=F(p,\{T_{2l,d}\})$ and this quantity
 is the exponential growth rate of the number of monomer-dimer tilings of
 the $d$-dimensional cubic grid having sides of length $n$, as $n$
 tends to infinity and the dimer density, (fraction of the maximum
 possible number $\frac{1}{2}n^d$ of dimers), in these tilings
 converges to a fixed number $p \in [0,1]$.  (See also \cite{FP1}.)
 We denote by this
 exponential growth by
 $h_d(p)$, and call it the $d$-dimensional monomer-dimer entropy
 of dimer density $p \in [0,1]$
 in the lattice $\Z^d$, see \cite{Ha1} and \cite{FP}.
 For $d=1$ the rate is
 known explicitly as a function $p$, whereas for $d > 1$ the exact rate
 is unknown.

 $h_d(p)$ can be estimated if one can estimate
 the quantities $f(p,\{G_l\})\le F(p,\{G_l\})$
 from below and above for any sequence $G_l\in
 \Gamma(n_l,r),l\in\N$.  Lower and upper estimates of
 $f(p,\{G_l\})\le F(p,\{G_l\})$ are conjectured in \cite{FKLM}
 and called the \emph{asymptotic lower matching conjecture}
 and \emph{asymptotic upper matching conjecture}, abbreviated
 here by ALMC and AUMC respectively.  For $r=2$ ALMC and AUMC
 are proved in \cite{FKM}.

 In this paper we apply our lower bounds on the sum of subpermanents
 of doubly stochastic matrices with $r$ nonzero
 elements in each row to obtain lower bound on $f(p,\{G_l\})$.
 For a fixed integer $r\ge 3$, we show
 the validity of ALMC for the densities $p_s=\frac{r}{r+s}$ where
 $s=0,1,\ldots$.  These inequalities yield new lower bounds
 for the $d$-dimensional monomer-dimer entropy of dimer density
 $h_d(p),p \in [0,1]$ in the lattice $\Z^d$.  In particular we
 obtain the best known lower bound for the three dimensional monomer dimer entropy
 $h_3$, which combined with the known upper bound in \cite{FP}
 gives the tight result $h_3\in [.7845, .7863]$.

 Next we discuss briefly the
 sum of $m\times m$ subhafnians of $2n\times 2n$ symmetric $B$ with nonnegative
 entries, denoted by $\haff_m B$.  For $0-1$ matrix $B$ this is equivalent to the number
 of $m$-matchings in a general graph on $2n$ vertices.
 We give Ryser-type algorithm for $\per_m A$ and $\haff_m B$.
 Unfortunately, the generating function $\x\trans B\x$,
 (the quadratic for associated with $B$), is positive hyperbolic
 if and only if the second eigenvalue of $B$ is nonpositive.
 We show that for any graph $G$, $\x\trans B(G)\x$ is positive
 hyperbolic if and only if $G$ is a complete $k$-partite graph.
 The last section is devoted to algorithmic applications related to
 $\per_m A$ and $\haff_m B$.

 We now list briefly the contents of this paper.
 In \S2 we discuss briefly the notion of positive hyperbolic polynomials
 and their properties that needed here.
 In \S3 we bring the generalized version of FT inequality for
 positive hyperbolic polynomials.  In \S4 we state and discuss the
 ALMC and AUMC.  In \S5 we give lower
 bounds on $f(p,\{G_l\})$.  We apply these bounds to verify the ALMC
 for a countable values of densities for each $r\ge 3$ as explained
 above.  In \S6 we discuss the notion of $\haff_m B$ and its
 connection to the quadratic form $\x\trans B\x$.
 In \S7 we discuss the Ryser-type formulas for $\perm_m A$ and $\haff_m B$.

 We thank Uri Peled for supplying us with the Figures 1 and 2.

 \section{Positive hyperbolic polynomials}

 \textbf{Definitions and Notations}

 \begin{enumerate}

 \item
 A vector $\x:=(x_1,\ldots,x_n)\trans\in \R^n$ is called
 \emph{positive} or \emph{nonnegative}, and denoted by $\x>\0:=
 (0,\ldots,0)\trans$ or $\x\ge \0$
 if $x_i>0$ or $x_i\ge 0$ for $i=1,\ldots,n$ respectively.
 A nonnegative vector $\x\ne \0$ is denoted by $\x\gvertneqq \0$.
 $\y\ge \x\iff \y-\x\ge \0$.
 The cone of all nonnegative vectors in $\R^n$ is denoted by
 $\R_+^n$.

 \item
 A polynomial $p=p(\x)=p(x_1,\ldots,x_n):\R^n
 \to\R$ is called \emph{positive hyperbolic} if the following conditions hold:

 \begin{itemize}
 \item $p$ is a homogeneous polynomial of degree $m\ge 0$.

 \item $p(\x)>0$ for all $\x>0$.

 \item $\phi(t):=p(\x+t\u)$, for $t\in \R$, has $m$-real $t$-roots for each
 $\u>\0$ and each $\x$.

 \end{itemize}

 \item For any polynomial $p:\R^n \to \R$ and any $\0\ne
 \u=(u_1,\ldots,u_n)\trans\in\R^n$ let $p_{\u}=p_{\u}(\x): \sum_{i=1}^n u_i\frac{\partial
 p}{\partial x_i}(\x)$.

 \item Let $\e_i:=(\delta_{i1},\ldots,\delta_{in})\trans \in \R^n,
 \;i=1,\ldots,n$ be the standard basis in $\R^n$.

 \item Let $\1:=(1,\ldots,1)\trans \in \R^n$ and denote by $J_n\in
 \R^{n\times n}$ the $n\times n$ matrix whose all entries are
 equal to $\frac{1}{n}$.
 \end{enumerate}

 The following lemma summarizes the basic properties of positive
 hyperbolic polynomials that needed here.
 \begin{lemma}\label{propphp}
 Let $p:\R^n \to \R$ be
 a positive hyperbolic polynomial of degree $m\ge 1$.
 Then the following properties hold:
 \begin{enumerate}
 \item
 Let $\u\gvertneqq \0, \x$ be fixed and denote
 $\phi(t)=p(\x+t\u)$.  Assume that $p(\u) >0$.  Then $\phi(t)$
 has $m$ real $t$ roots.  Furthermore
 $p_{\u}(\x)$ is a positive hyperbolic polynomial of degree
 $m-1$.
 $\y\ge \x\ge \0 \Rightarrow p(\y)\ge p(\x)\ge 0$.

 \item Let $\u\gvertneqq \0, \x\gvertneqq\0$ and
 assume that $p(\u)=0$.
 Then either $\phi(t)>0$ for all $t\ge 0$ or
 $p(\x)=0$ and $\phi(t)\equiv 0$.  Assume that $p(\x)>0$ and $\phi(t)$ is not a constant
 polynomial.
 Then all its roots are real and negative.
 If $p_{\u}$ is not a zero polynomial then $p_{\u}$
 is a positive hyperbolic of degree $m-1$.

 \item If $q((x_1,\ldots,x_{n-1})):=p((x_1,x_2,\ldots,x_{n-1},0))$ is not identically zero
 then $q$ is a positive hyperbolic of degree $m$ in $\R^{n-1}$.
 In particular,
 $r((x_1,\ldots,x_{n-1})):=$

 \noindent
 $\frac{\partial p}{\partial x_n}((x_1,\ldots,x_{n-1},0))$
 is either a zero polynomial or a positive hyperbolic polynomial in
 $n-1$ variables of degree $m-1$.

 \end{enumerate}

 \end{lemma}

 \proof

 We assume that $n>1$ otherwise all the results are trivial

 \begin{enumerate}

 \item The facts $\phi(t)$ has $m$ negative roots and
 $p_{\u}$ is positive hyperbolic is in \cite{Gar} or \cite{Kho}.
 Hence $p_{\u}\ge 0$ on $\R_+^n$.
 Let $\v\ge \0$ and assume
 that $\u_l>\0,l=1,\ldots$ be a sequence
 of vectors converging to $\v$.  Then $p_{\u_l}\to p_{\v}$.
 Therefore $p_{\v}\ge 0$ on $\R^n_+$.  In particular $p_{\e_i}(\x)=\frac{\partial
 p}{\partial x_i}(\x)\ge 0$ n $\R_+^n$ for $i=1,\ldots,n$.
 Thus $p(\x)$ is a nondecreasing function in each variable $x_i$
 on $R_+^n$.  Hence $\y\ge \x\ge 0 \Rightarrow p(\y)\ge p(\x)\ge
 0$.

 \item
 Recall Brunn-Minkowski theorem that
 $f:=p^{\frac{1}{m}}$ is convex on all positive vectors in $\R^n$
 \cite[Thm 2, 4)]{Kho}.  Since $f$ is continuous on
 $\R_+^n$, it follows that
 $f$ is nonnegative and convex on $\R^n_+$.
 If $p(\x)=0$ it follows that $f(\y)=0$ on the interval joining
 $\x,\u$.  Hence $\phi(t)\equiv 0 \iff p(\x)=0$.
 Assume that $p(\x)>0$.  Then $\phi(0)>0$.
 If $\phi(t)$ is constant then $\phi(t)=\phi(0)>0$.
 Assume that $\phi(t)$ is a nonconstant polynomial.
 Let $\u_l>\0,l=1,\ldots$ be a sequence
 of vectors converging to $\u$.
 Then $\phi_l(t):= p(\x+t\u_l) \to
 \phi(t)$.  Each $\phi_l$ has $m$ real negative zeros.  Consider the
 complex projective space $\C\bP^{n-1}$ with the homogenous
 coordinates $\z=(z_1,\ldots,z_n)\trans$.  The $m$-roots of
 $\phi_l(t)$ correspond to the intersection points of the line $L_l:=\{\z:=s\x+t\u_l,\;
 s,t\in \C\}$
 with the hypersurface $H=\{z\in \C\bP^n :\;p(\z)=0\}$ in $\C\bP^{n-1}$ of degree $m$.
 All these points are real and are located in the affine part of
 of the line $L_l$, i.e. $s=1$ and $t<0$.
 As $L_l\to L:=\{\z:=s\x+t\u,\;
 s,t\in \C\}$ it follows that $L_l\cap H\to L\cap H$, counting
 with the multiplicities if $\phi(t)$ is not constant.
 Thus all roots of $\phi(t)$ are nonnegative.  Since $\phi(0)>0$
 it follows that all the roots of $\phi(t)$ are negative.
 Hence $\phi(t) >0$ for $t\ge 0$.

 Assume that $p_{\u}\not\equiv 0$.  Recall that $p_{\u_l}\to
 p_{\u}$.  According to part \emph{1}  $p_{\u_l}(\y)\ge \p_{\u_l}(\x)\ge 0$
 for each pair $\y\ge \x\ge \0$.  Hence $p_{\u}(\y)\ge p_{\u}(\x)\ge
 0$.  If $p_{\u}(\y)=0$ for some $\y>\0$ then for any $\x\ge 0$,
 there exists $t>0$ such that $\y\ge t\x\ge 0$.  Hence $t^{m-1}\p_{\u}(\x)=p_{\u}(t\x)=0
 \Rightarrow p_{\u}(\x)=0 \Rightarrow p_{\u}\equiv 0$ contrary to
 our assumption.  Thus $p(\y)>0$.  As $p_{\u_l}(\w+t\y) \to
 \psi(t):=p_{\u}(\w+ty)$ and $\psi(t)$ is a polynomial of degree
 $m-1$, it follows that $\psi(t)$ has $m-1$ real roots for each $\w\in\R^n$.
 Thus $p_{\u}$ is positive hyperbolic.

 \item
 Since $p(\y)\ge p(\x)\ge 0$ for $\y\ge \x$ on $\R_+^n$
 it follows that $q(\y_1)\ge q(\x_1)\ge 0$ on $\R_+^{n-1}$,
 where
 $\y=(\y_1\trans,0)\trans, \x=(\x_1\trans,0)\trans$.
 Since $q$ is a homogeneous of degree $m$, from the arguments in
 part \emph{2} it follows that either $q(\y_1)>0$ for each $\y_1>\0$
 or $q\equiv 0$.  Assume that $p(\y)=q(\y_1)>0$.  Then according
 to part \emph{1}, $p(\x+t\y)=q(\x_1+t\y_1)$ has $m$ real roots in $t$.
 Assume that $r\not\equiv 0$.  Hence $p_{\e_n}\not \equiv 0$.
 So $p_{\e_n}$ is positive hyperbolic by \emph{2}, and by previous
 arguments $r$ is positive hyperbolic on $\R^{n-1}$.

 \end{enumerate}
 \qed

 The following propositions are well known and we bring their proof
 for completeness.
 \begin{prop}\label{nonnegcphp}  Let $p:\R^n \to \R$ be
 a positive hyperbolic polynomial of degree $m\ge 1$.
 Then the coefficient of each monomial in $p$ is nonnegative.
 \end{prop}
 \proof
 Let $p$ be a positive hyperbolic polynomial of degree $m\ge 1$.
 Let $\u>\0$.  Part \emph{1} of the above lemma yields that $p_{\u}$ is positive
 hyperbolic of degree $m-1$.

 We prove the proposition by induction on $m$.
 Let $m=1$.  Assume that $\u>\0$.  Then $p_{\u}=\nabla p \u>0$.
 Hence $\nabla p \ge \0$ and the corollary holds.

 Assume that the result hold for $m=l\ge 1$.
 Let $p$ be a positive hyperbolic polynomial of degree $m=l+1$.
 Let $\e_i=(\delta_{i1},\ldots,\delta_{in})\trans\in\R^n$ for $i=1,\ldots,n$.
 Let $\u_j>\0, j=1,\ldots,$ and assume that
 $\lim_{j\to\infty}\u_j=\e_i$.
 Hence $p_{\u_j}$ is positive
 hyperbolic of degree $l$.  By the induction hypothesis the coefficients of all
 monomials of $p_{\u_j}$ are nonnegative.  Let $j\to\infty$ and deduce
 that the coefficients of all monomials of $p_{\e_i}$ are
 nonnegative.  Hence the coefficients of all monomials of $p$ which
 include the variable $x_i$ are nonnegative.  Let $i=1,\ldots,n$ to
 deduce that the coefficients of all monomials of degree one at
 least are nonnegative.  As $p(\0)\ge 0$ we deduce the proposition.
 \qed

 \begin{prop}\label{logconsumexp}  For $\x=(x_1,\ldots,x_n)\trans\in\R^n$
 let $e^{\x}:=(e^{x_1},\ldots,e^{x_n})\trans>\0$.  Let $p:\R^n \to \R$
 be a nonzero polynomial with such that the coefficient of
 each monomial is nonnegative.  Then $\log p(e^{\x})$ is a convex function
 on $\R^n$.  Let $L\subseteq \R^n$ be a line in $\R^n$.  Then the
 restriction of $\log p(e^{\x})$ to $L$ is either an
 affine function or a strictly convex function.
 \end{prop}

 \proof   The convexity of $\log p(e^{\x})$ can be found in \cite{Kin}.
 Note that $f:=p(e^{\x})|L$ is of the form in $\sum_{i=1}^k a_ie^{b_it}$, where each $a_i >0$
 and $b_1 > b_2>\ldots > b_k$.
 If $k=1$ then $\log f(t)$ is $\log a_1 + b_1t$.  Otherwise it
 straightforward to show that $\frac{f'}{f}$ is increasing on $\R$.
 \qed

 \begin{corol}\label{convlogpshyp}  Let $p:\R^n \to \R$ be
 a positive hyperbolic polynomial of degree $m\ge 1$.
 Then $\log p(e^{\x})$ is a convex function
 on $\R^n$.  Let $L\subseteq \R^n$ be a line in $\R^n$.  Then the
 restriction of $\log p(e^{\x})$ to $L$ is either an
 affine function or a strictly convex function.
 \end{corol}

 \textbf{Examples of positive hyperbolic polynomials}

 \begin{enumerate}
 \item
 Let $A=(a_{ij})_{i=j=1}^{m,n}\in \R^{m\times n}$
 be a nonnegative matrix, denoted by $A\ge \0$, where each
 row of $A$ is nonzero.  Fix a positive integer $k\in [1,m]$.
 Then
 \begin{equation}\label{exhpk}
 p_{k,A}(\x):=\sum_{1\le i_1<\ldots i_k\le m} \prod_{j=1}^k(A\x)_{i_j},
 \x\in\R^n,
 \end{equation}
 is positive hyperbolic of degree $k$ in $n$ variables.

 \item  Let $A_1,\ldots,A_n\in \C^{m\times m}$ hermitian,
 nonnegative definite matrices such that $A_1+\ldots+A_n$ is
 a positive definite matrix.  Let $p(\x)=\det \sum_{i=1}^n x_i
 A_i$.  Then $p(\x)$ is positive hyperbolic.

 \end{enumerate}

 \proof
 \begin{enumerate}
 \item
 First note that $p_{k,A}(\x)>0$ for $\x>\0$.
 The hyperbolicity of $p_{m,A}$ and $p_{1,A}$ is obvious.
 Assume that $k\in (1,m)$.  Let
 $\z=(z_1,\ldots,z_{n+m-k})\trans\in\R^{n+m-k}$ and define
 $P(\z):=\prod_{i=1}^m (\sum_{j=1}^n a_{ij}z_j+\sum_{j=n+1}^{n+m-k}
 z_j)$.  Then
 $$p_{k,A}(\x)={m \choose k}^{-1}\frac{\partial^{m-k} P}{\partial z_{n+1}\ldots
 \partial z_{n+m-k}}((x_1,\ldots,x_n,0,\ldots,0)).$$
 Hence by Lemma \ref{propphp} $p_{k,A}$ positive hyperbolic.

 \item  This is a standard example and the proof is
 straightforward.

 \end{enumerate}
 \qed

 Let $p(\x):\R^n\to \R$ be a positive hyperbolic polynomial of degree $m\ge 1$.
 As in \cite{Gur1} define
 \begin{equation}\label{defcap}
 \Cp p:=\inf_{\x>\0, x_1\ldots x_n=1} p(\x)=\inf_{\x>\0}
 \frac{p(\x)} {(x_1\ldots x_n)^{\frac{m}{n}}}.
 \end{equation}
 It is possible that $\Cp p=0$.  For example let
 $p=x_1^{m_1}\ldots x_n^{m_n}$ where $m_1,\ldots,m_n$ are
 nonnegative integer whose sum is $m$ and $(m_1,\ldots,m_n)\ne
 k\1$.

 \begin{prop}\label{excappos}
 Let $A\in\R^{n\times n}$ be a doubly stochastic matrix.
 Let $p_{k,A},k\in [1,n]$ be positive hyperbolic defined as
 part 1 of the above example.  Then $\Cp p_{k,A}={n \choose k}$.
 Let $B\in\R^n$ be a matrix with positive entries.
 Then there exists two positive definite diagonal matrices
 $D_1,D_2$, unique up to $tD_1,t^{-1}D_2, t>0$, such that $A:=D_1BD_2$ is a doubly
 stochastic matrix \cite{Sin}.  Let $p_{n,B}$ be defined as above.
 Then $\Cp p_{n,B}=\frac{1}{\det D_1D_2}$.
 \end{prop}
 \proof  Consider first $p_{n,A}$.  Since $A$ is row stochastic
 $p_{n,A}(\1)=1$.  Hence $\Cp p_{n,A}\le 1$.  Let
 $\u=(u_1,\ldots,u_n)\trans \gvertneqq \0$ be a probability
 vector.  Then for any $\x=(x_1,\ldots,x_n)>\0$ the generalized
 arithmetic-geometric inequality states $\u\trans \x\ge
 \prod_{i=1}^n x_i^{u_i}$.  Use this inequality for each
 $(A\x)_i$.  The assumption that $A$ is doubly stochastic
 yields that $p_{n,A}\ge x_1\ldots x_n \Rightarrow \Cp p_{n,A}\ge
 1$.  Hence $\Cp  p_{n,A}=1$.

 Let $k\in [1,n)$.  Then $p_{k,A}(\1)={n \choose k}$.  Hence
  Hence $\Cp p_{k,A}\le {n \choose k}$.
 Apply the arithmetic-geometric inequality to
 $\frac{p_{k,A}}{{n \choose k}}$ to deduce that
 $p_{k,A}\ge {n \choose k} p_{n,A}^{\frac{m}{n}}$.
 Hence $\Cp p_{k,A}\ge {n \choose k}$.

 It is straightforward to show that  $\frac{p_{n,B}(\x)}{x_1\ldots
 x_n}=\frac{p_{n,A}(\y)}{\det(D_1D_2) y_1\ldots y_n}$, where
 $\y=D_2^{-1}\x$.  Hence $\Cp p_{n,B}=\frac{1}{\det D_1D_2}$.
 \qed

 The following result is taken from \cite{Gur1}.

 \begin{lemma}\label{derest}
 Let $k\ge 1$ be an integer, $\u:=(u_1,\ldots,u_k)\trans> \0, \v:=(v_1,\ldots,v_k)\trans >\0$
 and define $f(t):=\prod_{i=1}^k (u_it+v_i)$.  Let
 $K(f):=\inf_{t>0} \frac{f(t)}{t}$.  Then
 $f'(0)=K$ for $k=1$ and $f'(0)\ge (\frac{k-1}{k})^{k-1} K$ for $k\ge 2$.
 For $k\ge 2$ equality holds if and only if
 $\frac{v_1}{u_1}=\ldots=\frac{v_k}{u_k}$.
 \end{lemma}

 \proof
 We can assume WLOG that $f(0)=1$ ; i.e. that
 $f(t):=\prod_{i=1}^k (a_{i}t + 1) , a_i = \frac{u_{i}}{v_{i}}$.
 Using the arithmetic-geometric means inequality we get that
 $$
 Kt \leq f(t) \leq p(t) =: (1 + \frac{f'(0)}{k} t)^k , t \geq 0.
 $$
 Therefore , by doing basic calculus ,
 $$
 K \leq \inf_{t>0} \frac{p(t)}{t} = f'(0)(\frac{k}{k-1})^{k-1} ,
 $$
 which finally gives the desired inequality
 $$
 f'(0)\ge (\frac{k-1}{k})^{k-1} K , k \geq 2.
 $$

 It follows again from arithmetic-geometric means inequality that
 the equality holds if and only if $a_1 = \frac{v_1}{u_1}=\ldots= a_k = \frac{v_k}{u_k}$.
 \qed

 \textbf{Definition.}
 Let $p:\R^n\to \R$ be a positive hyperbolic polynomial of degree
 $m$.  For each integer $i\in [0,n]$ the \emph{i-th degree} of $p$
 is the integer $r_i\in [0,m]$ such that
 $$\frac{\partial^{r_i} p}{\partial
 x_i^{r_i}}(x_1,\ldots,x_{i-1},0,x_{i+1},\ldots,x_n)\not
 \equiv 0, \textrm{ and }
 \frac{\partial^{r_i+1} p}{\partial
 x_i^{r_i+1}}(x_1,\ldots,x_{i-1},0,x_{i+1},\ldots,x_n)\equiv 0.$$
 Let $\deg_i p:=r_i$ for $i=1,\ldots,n$.

 The following proposition follows straightforward from part \emph{3}
 of Lemma \ref{propphp}.
 \begin{prop}\label{ranklem}  Let $p:\R^n \to \R$ be a positive
 hyperbolic polynomial of degree $m$.  Let $i\in [1,n]$ be an integer.  Then
 \begin{enumerate}
 \item $\deg_i p= 0 \iff
 p(\x)=(p(x_1,\ldots,x_{i-1},0,x_{i+1},\ldots,x_n))$.
 \item For each integer $j\in [0,\deg_i p]$ $\frac{\partial^{j} p}{\partial
 x_i^{j}}(x_1,\ldots,x_{i-1},0,x_{i+1},\ldots,x_n)$ is a positive
 hyperbolic polynomial of degree $m-j$.
 \item For each integer $j\in [1,n], j\ne i$,
 $$\deg_j
 \frac{\partial p}{\partial
 x_i}(x_1,\ldots,x_{i-1},0,x_{i+1},\ldots,x_n) \le \min (\deg_j p,
 n-1).$$

 \end{enumerate}
 \end{prop}

 The following result is crucial for the proof of
 a generalized Friedland-Tverberg inequality and is due essentially to the second author in \cite{Gur1}.

 \begin{lemma}\label{capest}  Let $p:\R^n\to \R$ be a positive hyperbolic
 polynomial of degree $m\ge 1$.  Assume that $\Cp p>0$.
 Then $\deg_i p \ge 1$ for $i=1,\ldots, n$.
 For $m=n\ge 2$
 $$\Cp \frac{\partial p}{\partial
 x_i}(x_1,\ldots,x_{i-1},0,x_{i+1},\ldots,x_n)\ge (\frac{\deg_i p
 -1}{\deg_i p})^{\deg_i p-1} \Cp p \textrm{ for }
 i=1,\ldots,n,$$
 where $0^0=1$.
 \end{lemma}

 \proof  It is enough to prove the result for $i=n$.  Suppose to
 the contrary that $p$ does not depend on $x_n$.  Then let
 $\x(t)=(1,\ldots,1,t)\trans$ and $t\to\infty$ in (\ref{defcap}) to deduce
 that $\Cp p=0$ contrary to our assumption.

 Assume that $m=n>1$.  Let $k=\deg_n p\ge 1$.
 Let $\x_0:=(x_1,\ldots,x_{n-1},0)\trans, \x_1:=(x_1,\ldots,x_{n-1})\trans$.
 Proposition \ref{ranklem}
 yields that
 $g(\x_1):=\frac{\partial^{k} p}{\partial
 x_i^{k}}(\x_0)$ is a positive hyperbolic function
 in $n-1$ variables of degree $m-l$.  Hence $g(\x_1)>0$ for
 $\x_1>\0$.  Thus for $\x_1>\0$
 \begin{equation}\label{capest1}
 p(\x_0+t\e_n)=k!g(\x_1)t^l+\ldots=k!g(\x_1)\prod_{i=1}^k
 (t+\lambda_i(\x_1)), \quad \lambda_i(x)>0, \textrm{ for }i=1,\ldots,k.
 \end{equation}
 The second equality follows from part \emph{2} of Lemma
 \ref{propphp}.  Assume in addition that $x_1\ldots x_{n-1}=1$.
 Then $\inf_{t>0} \frac{p(\x_0+t\e_n)}{t}\ge \Cp p$.  Apply
 Lemma \ref{derest} to the right-hand side of (\ref{capest1}) to deduce that
 $\frac{\partial p}{\partial
 x_n}(\x_0)\ge (\frac{k-1}{k})^{k-1} \Cp p$.
 Since we assumed that $x_1\ldots x_{n-1}=1$ it follows that
 $\Cp\frac{\partial p}{\partial
 x_n}(\x_0)\ge (\frac{k-1}{k})^{k-1} \Cp p$.
 \qed

 \begin{rem}\label{concave}
 Lemma \ref{derest} , which is simple but crucial , is a particular case of the
 following general result :\\
 {\it Let $f : [0, \infty) \rightarrow R_{+}$ be a nonnegative function differentiable at zero
 from the right ; $K = \inf_{t>0} \frac{f(t)}{t}$. If $k \geq 1$ and $f^{\frac{1}{k}}$ is concave then
 $f'(0) \geq (\frac{k-1}{k})^{k-1} K$. On the other hand if $k \geq 1$ and $f^{\frac{1}{k}}$ is convex
 then $f'(0) \leq (\frac{k-1}{k})^{k-1}K$.}

 \end{rem}
 \section{Friedland-Tverberg inequality}

 \begin{theo}\label{frtver}  Let $p:R^n \to \R$ be  positive
 hyperbolic of degree $m\in [1,n]$.
 Assume that $\deg_i p\le r_i \in [1,m]$ for $i=1,\ldots, n$.
 Rearrange the sequence $r_1,\ldots,r_n$ in an increasing
 order $1\le r_1^*\le r_2^*\le \ldots \le r_n^*$.  Let $k\in [1,n]$ be
 the smallest integer such that $r_k^*> m-k$.
 Then

 \begin{eqnarray}
 \sum_{1\le i_1 <\ldots<i_m\le n} \frac{\partial ^m p}{\partial
 x_{i_1}\ldots\partial x_{i_m}} (\0) \ge
 \nonumber\\
 \frac{n^{n-m}}{(n-m)!}\frac{(n-k+1)!}{(n-k+1)^{n-k+1}}\prod_{j=1}^{k-1}
 (\frac{r_j^*+n-m-1}{r_j^*+n-m})^{r_j^*+n-m-1} \Cp p.
 \label{frtver1}
 \end{eqnarray}
 (Here $0^0=1$, and the empty product for $k=1$ is assumed to be
 $1$.)
 If $\Cp>0$ and $r_i=m$ for $i=1,\ldots,m$ equality holds if and only if
 $p=C(\frac{x_1+\ldots+x_n}{n})^m$ for each $C>0$.
 \end{theo}
 \proof  Suppose that $\Cp p=0$.  Then part \emph{3} of Lemma
 \ref{propphp} yields that the left-hand side of (\ref{frtver1})
 is nonnegative and the theorem holds in this case.

 Clearly, it is enough to assume the
 case $\Cp p=1$.  The case $m=n$ is essentially proven  in \cite{Gur1} and we
 repeat its proof for the convenience of the reader.
 Permute the coordinates of $x_1,\ldots,x_n$ such that $\deg_n
 p=\min_{i\in[1,n]}\deg_i p\le r_1^*$.
 Assume that $\deg_n p=l$.  Then Lemma \ref{capest} yields that
 $r((x_1,\ldots,x_{n-1})):=\frac{\partial p}{\partial x_n}((x_1,\ldots,x_{n-1},0))$
 is positive hyperbolic of degree $n-1$ and $\Cp r\ge (\frac{l-1}{l})^{l-1} \Cp p$.
 Since the sequence $(\frac{i-1}{i})^{i-1}, i=1,\ldots,$ is
 decreasing to have the lowest possible lower bound we have to assume
 $l=r_1^*$.  Suppose first that $r_1^*=n$.
 Repeating this process $n$ times we get that
 $$\frac{\partial ^n p}{\partial
 x_{1}\ldots\partial x_{n}} (\0)\ge \Cp p\prod_{j=2}^n
 (\frac{j-1}{j})^{j-1}=\frac{n!}{n^n} \Cp p.$$
 This inequality to corresponds to the case $r_i^*=n$ for
 $i=1,\ldots,n$.
 The equality case is discussed in \cite{Gur1}.

 Let $m\in [1,n-1]$.  Put $P(\x)=p(\x)(\frac{1}{n}\sum_{i=1}^n
 x_i)^{n-m}$.  Clearly, $P$ is positive hyperbolic of degree $n$.
 Since $\frac{1}{n}\sum_{i=1}^n x_i\ge
 (x_1\ldots x_n)^{\frac{1}{n}}$ for each $\x\ge 0$, it follows that $\Cp P\ge \Cp p$.
 Apply (\ref{frtver1}) to $P$ for $m=n$ to deduce (\ref{frtver1})
 in the general case.  Since the equality case for $P$ holds
 if and only $P=(\frac{1}{n}\sum_{i=1}^n x_i)^{n}$ it follows that
 the equality in (\ref{frtver1}) holds if and only if $p=(\frac{1}{n}\sum_{i=1}^n
 x_i)^{m}$.  \qed

 Let $A\in\R^{n\times n}$ be a doubly stochastic matrix.
 Apply this theorem to $p_{m,A}$ defined Proposition \ref{excappos}
 to deduce the Friedland-Tverberg inequality for the sum of all
 $m\times m$ permanents of $A$:
 \begin{corol}\label{frteverds}  Let $A\in\R^{n\times n}_+$ be a doubly stochastic
 matrix.  Then $\perm_m A\ge {n \choose m}^2 \frac{m!}{n^m}$ for any $m\in [2,n]$.
 equality holds if and only $A=J_n$.
 \end{corol}

 \begin{theo}\label{gurschr}  (\textrm{Gurvits})  Let  $A\in\R^{n\times n}_+$ be a doubly stochastic
 matrix, such that each column contains at most $r\in [1,n]$ nonzero
 entries.  Then
 \begin{equation}\label{gurschr1}
 \per A \ge \frac{r!}{r^{r}}\big(\frac{r-1}{r}\big)^{(r-1)(n-r)}=\frac{r!}{r^r}\big(\frac{r}{r-1}
 \big)^{r(r-1)}\big(\frac{r-1}{r}\big)^{(r-1)n}.
 \end{equation}
 \end{theo}
 \proof  Note that for $p(\x)=\prod_{i=1}^n (A\x)_i$ we have that
 $\deg_i p=r$ for $i=1,\ldots,n$.  Apply (\ref{frtver1}) to this
 case, i.e. $m=n$, $r_j^*=r, j=1,\ldots,n$ and $k=n-r+1$ to deduce
 the theorem.
 \qed

 \section{The ALMC and AUMC}

 Let $G=(V,E)$ be a general graph with the set of vertices $V$ and
 edges $E$.
 A \emph{matching} in $G$ is a
 subset $M \subseteq E$ such that no two edges in $M$ share a common
 endpoint. The endpoints of the edges in $M$ are said to be
 \emph{covered} by $M$. We can think of each edge $e = (u,v) \in M$
 as occupied by a \emph{dimer}, consisting of two neighboring atoms
 at $u$ and $v$ forming a bond, and of each vertex not covered
 by $M$ as a \emph{monomer}, which is an atom not forming any bond.
 For this reason a matching in $G$ is also called a
 \emph{monomer-dimer cover} of $G$. If there are no monomers, $M$ is
 said to be a \emph{perfect matching}.  Note that if a perfect
 matching exists then $\#V$ is even.
 A matching $M$ with
 $\# M = k$ is called an \emph{$k$-matching}. We denote by
 $\phi_G(k)$ be the number of $k$-matchings in $G$ (in particular
 $\phi_G(0) = 1$), and by $\Phi_G(x):= \sum_k \phi_G(k) x^k$ the
 matching generating polynomial of $G$.  It is known that all the
 roots of matching polynomial are real negative numbers \cite{LP}.

 Let $G$ be a bipartite graph $G=(V,E)$, where $V=V_1\cup V_2$ is
 the set of vertices of $G$ and $E$ is the set of edges that
 connect vertices in $V_1$ to vertices in $V_2$.  Assume that
 $\#V_1=\#V_2=n$. We identify $V_1$ and $V_2$ with
 $\an{n}:=\{1,\ldots,n\}$, where the vertices in $V_1$ and $V_2$ are colored
 with colors black and white respectively.
 Then $G$ is represented by $0-1$ $n\times n$ matrix
 $A(G)=A=(a_{ij})\in \{0,1\}^{n\times n}$, where $a_{ij}=1$
 if and only if the black vertex $i$ is connected to the white
 vertex $j$.  It is convenient to consider multi bipartite graphs.
 Thus, the entries of the representation matrix $A(G)=(a_{ij})\in
 \Z_+^{n\times n}$ are nonnegative integers, where $a_{ij}$ is the
 number of edges from the black vertex $i$ to the white vertex $j$.

 It is straightforward to show that
 \begin{equation}\label{matperid}
 \phi_G(k)=\perm_k A(G), \quad k=0,\ldots,n, \textrm{ where }
 \perm_0 A:=1 \textrm{ for any } A\in \R^{n\times n}.
 \end{equation}

 Let $\Gamma(n,r)$ be the set of bipartite $r$-regular multi-graphs,
 (each vertex has degree $r$), with $n:=\frac{\#V}{2}$.
 Let $\Delta(n,r)$ be the set
 by an $n\times n$ nonnegative matrices $A$ with integer
 entries, such that the sum of each row and column is $r$.
 Then each $G\in \Gamma(n,r)$ is represented by $A\in
 \Delta(n,r)$ and vice versa.
 Note for each $A\in \Delta(n,r)$
 the matrix $\frac{1}{r}A$ is doubly stochastic.
 Corollary \ref{frteverds} yields:
 \begin{equation}\label{frtverbinmat}
 \phi_G(m)\ge {n\choose m}^2\frac{m!r^m}{n^m} \quad \textrm{ for any }
 G\in  \Gamma(n,r).
 \end{equation}

 Note that the symmetric group $S_n$ on $n$ elements, presented
 as the group of permutation matrices $\Pi_n\subset \{0,1\}^{n\times
 n}$ acts from the left and from the right on $\Delta(n,r)$, i.e.
 $P\Delta(n,r)=\Delta(n,r)P$ for each $P\in\Pi_n$.
 These actions are equivalent to the action of $S_n$ on $V_1$ and
 $V_2$ respectively.

 There is a standard probabilistic model on $\Gamma(n,r)$,
 which assigns a fairly natural probability measure $\nu(n,r)$ on
 $\Gamma(n,r)$ \cite{LP}.  The measure $\nu(n,r)$
 is invariant under the action of $S_n$ on $V_1$ and $V_2$
 as explained above.  By abuse of the notation we view $\nu(n,r)$
 also a probability measure on $\Delta(n,r)$, which is invariant under the left and the
 right action of $\Pi_n$.
 The following result is proven in \cite{FKM}:

 \begin{theo}\label{mexpcval}  Let $\nu(n,r)$ be the probability
 measure defined above.  Then
 \begin{equation}\label{mexpcval1}
 E_{\nu(n,r)}(\phi(G,m))=  E_{\nu(n,r)}(\perm _m A)=\frac{{n \choose m}^2
 r^{2m}m!(rn-m)!}{(rn)!}, \quad m=0,\ldots,n.
 \end{equation}

 In particular, let $k_n \in [0,n]$, $n = 1,2,\ldots$
      be a sequence of integers with $lim_{n \to \infty}
      \frac{k_n}{n} = p\in (0,1]$.  Then
  \begin{equation}\label{lasavlim}
       \lim_{n \to \infty} \frac{\log E_{\nu(n,r)}(\phi(G,k_n))}{2n}  =   gh_r(p),
      \end{equation}
  where
 \begin{equation}\label{defghrp}
      gh_r(p) := \textstyle{\frac{1}{2}}\left(p \log r -p\log p - 2(1-p)\log (1-p) +(r-p)\log \left(1
      -\frac{p}{r}\right)\right),
 \end{equation}

 \end{theo}
 The case $m=n$ in (\ref{mexpcval1}) is given in \cite{LP}.

 Fix as subset $J\subset \an{n}$ of cardinality $m$: $\#J=m$.
 For $G\in\Gamma(n,r)$ let $\phi(G,J)$ be all $m$-matching in $G$
 that cover the set $J\subset V_2$.  For $A=(a_{ij})\in \R^{n\times n}$ and $I\subset \an{n}$ let
 $A[J|I]$ be the submatrix $(a_{ij})_{i\in I,j\in J}$.  Denote
 $$\perm_m A[\an{n}|J]=\sum_{I\subset \an{n},\#I=m}\perm _m
 A[I|J].$$
 Then $\phi(G,J)=\perm _m A(G)[\an{n}|J]$.

 Use the invariance of $\nu(n,r)$ under the action
 of $S_n$ on $V_2$ and the fact that there are
 $n\choose r$ distinct subsets $J\subset \an{n}$ of cardinality $m$
 to obtain:
 \begin{corol}\label{mexpcvalcor}  Let $\nu(n,r)$ be the probability
 measure defined above.  Then for any $J\subset \an{n}, \#J=m$
 \begin{equation}\label{mexpcval1cor}
 E_{\nu(n,r)}(\phi(G,J))=  E_{\nu(n,r)}(\perm _m A[\an{n}|J])=\frac{{n \choose m}
 r^{2m}m!(rn-m)!}{(rn)!}, \quad m=0,\ldots,n.
 \end{equation}

 \end{corol}

 The following conjecture is stated in \cite{FKLM}.

  \begin{con}[The Asymptotic Lower Matching Conjecture]
      \label{lasmc} \mbox{}\newline
      For $r \geq 2$, let $G_n=(V_n,E_n)$, $n=1,2,\ldots$ be a
      sequence of finite $r$-regular bipartite graphs with $\#V_n \to
      \infty$. Let $k_n \in [0,\frac{\#V_n}{2}]$, $n = 1,2,\ldots$
      be a sequence of integers with $lim_{n \to \infty}
      \frac{2k_n}{\#V_n} = p\in (0,1]$.  Then
      \begin{equation}\label{lasmcp}
       \liminf_{n \to \infty} \frac{\log \phi_{G_n}(k_n)}{\#V_n} \geq gh_r(p).
      \end{equation}
 \end{con}
 For $r=1$ this conjecture holds trivially.  For $r=2$ this
 conjecture is proved in \cite{FKM}.
 The inequality (\ref{frtverbinmat}) implies that under the
 conditions of Conjecture \ref{lasmc} the following inequality
 holds, see \cite{FP}
 \begin{equation}\label{lasmbdfp}
       \liminf_{n \to \infty} \frac{\log \phi_{G_n}(k_n)}{\#V_n} \geq
       fh_r(p),
      \end{equation}
 where
 \begin{equation}\label{FTLowerBound}
  fh_r(p) := {\textstyle \frac{1}{2}} ( -p\log p - 2(1-p)\log(1-p) +
  p\log r -p ).
 \end{equation}

  As usual, we denote by $\R[x]$ the algebra of polynomials in $x$
  with real coefficients, by $0 \in \R[x]$ the zero polynomial, and
  by $\R_+[x] \subset \R[x]$ the subalgebra of polynomials with
  non-negative coefficients. We partially order $\R[x]$ by writing,
  for $f, g \in \R[x]$, $g \succeq f$ when $g - f \in \R_+[x]$, and
  $g \succ f$ when $g - f \in \R_+[x] \setminus \set{0}$. Clearly,
  if $g_1 \succeq f_1 \succ 0$ and $g_2 \succeq f_2 \succ 0$, then
  $g_1g_2 \succ f_1 f_2$ unless $g_1=f_1$ and $g_2=f_2$.

 Let $qK_{r,r}$ denote the union of $q$ complete bipartite graphs $K_{r,r}$
 having $r$ vertices of each color class.
 It is straightforward to show that any finite graphs $G,G'$ satisfy
  \begin{equation}\label{prodgenmatp}
   \Phi_{G \cup G'}(x)  =\Phi_G(x)\Phi_{G'}(x),
  \end{equation}
  and that
  \begin{equation}\label{matkrr}
       \Phi_{K_{r,r}}(x) = \sum_{k=0}^r \binom{r}{k}^2 k!\,x^k.
  \end{equation}
  The following conjecture is stated in \cite{FKLM}
  \begin{con}[The Upper Matching Conjecture]\label{upmatcon}
       Let $G$ be a bipartite $r$-regular graph on $2qr$ vertices
       where $q,r \geq 2$. Then  $\Phi_G \preceq \Phi_{q K_{r,r}}$,
       equality holding only if $G = q K_{r,r}$.
  \end{con}
 For $k=2$ this conjecture is proved in \cite{FKM}.
 The above conjecture implies the following Asymptotic Upper
 Matching Conjecture \cite{FKLM}. Denote by $K(r)$ be the countably
 infinite union of $K_{r,r}$.
 Let $ P_{K(r)}(t)$ and  $h_{K(r)}(p), p\in [0,1]$ be the pressure and the
 $p$-matching entropy associated
 and the with $K(r)$ \cite{FP1}:
  \begin{equation}\label{presK}
       P_{K(r)}(t) = \frac{\log \sum_{k=0}^r \binom{r}{k}^2 k!\, e^{2kt}}{2r},
       \qquad t \in \R.
  \end{equation}

  \begin{equation}\label{hKp}
       h_{K(r)}(p(t)) = P_{K(r)}(t)
       - t p(t), \qquad t \in \R
  \end{equation}
  where
  \begin{equation}\label{pt}
         p(t) = P'_{K(r)}(t) = \frac{\sum_{k=0}^r \binom{r}{k}^2 k!\, (2k)e^{2kt}}
       {2r\sum_{k=0}^r \binom{r}{k}^2 k!\, e^{2kt}}, \qquad t \in \R.
  \end{equation}
 \begin{con}[The Asymptotic Upper Matching Conjecture]\label{upasmatcon}
  \mbox{}\\
 For $r \geq 2$, let $G_n=(V_n,E_n)$, $n=1,2,\ldots$ be a
      sequence of finite $r$-regular bipartite graphs with $\#V_n \to
      \infty$. Let $k_n \in [0,\frac{\#V_n}{2}]$, $n = 1,2,\ldots$
      be a sequence of integers with $lim_{n \to \infty}
      \frac{2k_n}{\#V_n} = p\in (0,1]$.  Then
      \begin{equation}\label{uasmcp}
       \limsup_{n \to \infty} \frac{\log \phi_{G_n}(k_n)}{\#V_n} \leq h_{K(r)}(p).
      \end{equation}
      Equality case holds for the sequence $qK_{r,r},
      q=1,2,\ldots$.
 \end{con}
 For $r=2$ the AUMC is proven in \cite{FKM}.
 For $p=1$ and any $r\in\N$ the AUMC follows from the proof of Minc
 conjecture by Bregman \cite{Bre}.  Some computations performed in
 \cite{FKLM} support the ALMC and AUMC.

\begin{figure}[h]
    \begin{center}
    \includegraphics[width=0.75\textwidth]{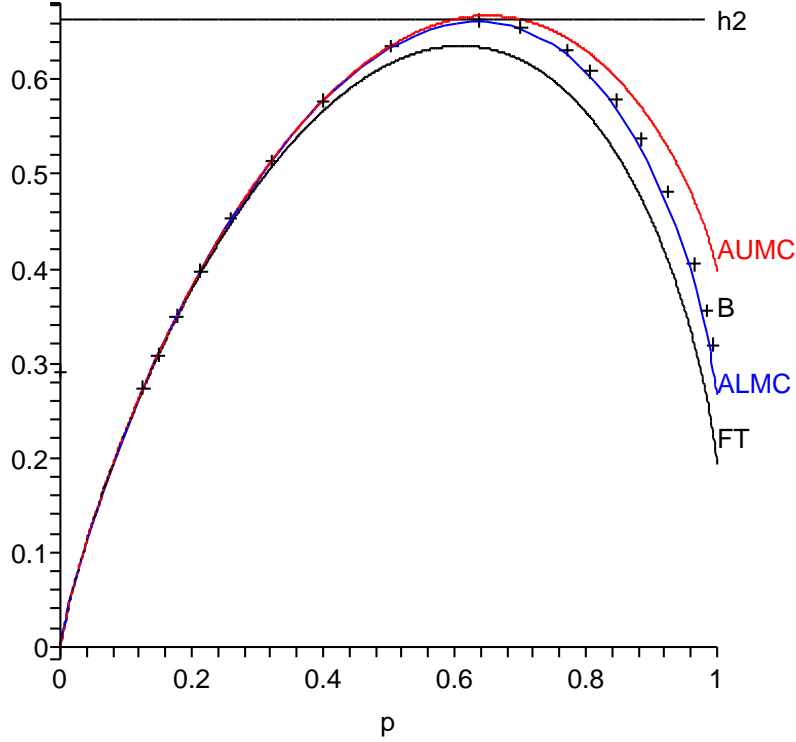}
    \caption{Monomer-dimer tiling of the $2$-dimensional grid:
    entropy as a function of dimer density. FT is the
    Friedland-Tverberg lower bound, h2 is the true monomer-dimer
    entropy. B are Baxter's computed values. ALMC is the Asymptotic
    Lower Matching Conjecture. AUMC is the entropy of a countable
    union of $K_{4,4}$, conjectured to be an upper bound by the
    Asymptotic Upper Matching Conjecture.}\label{Fig:h2bounds}
    \end{center}
  \end{figure}

 The following plots
 illustrating the Asymptotic Matching Conjectures for $r=4,6$.
 Let $C_n$ a cycle on $n$ points, and let
 $T_{n,d}=(V_n,E_n):=\underbrace{C_{n}\times\ldots
 \times C_n}_d, n=3,\ldots$ be a sequence of $d$ dimensional
 torii.  Note that each $T_{n,d}$ is $2d$ regular graph.  It is a classical result
 that the following limit exists for any $p\in[0,1]$:
 \begin{equation}\label{hdpdef}
      \lim_{n \to \infty} \frac{\log \phi_{T_{n,d}}(k_n)}{\#V_n}=h_d(p),
      \quad p\in [0,1].
     \end{equation}
 $h_d(p)$ is the $d$-dimensional monomer-dimer entropy of dimer density $p \in [0,1]$
 in the lattice $\Z^d$ \cite{Ha1} and \cite{FP}.
 In this case we use the notation $h_d := \max_{p\in[0,1]} h_d(p)$, (the quantities $h_d$
 and $\tilde h_d := h_d(1)$ are called the $d$-monomer-dimer entropy
 and the $2$-dimer entropy, respectively, in \cite{FP}).
 Figure~\ref{Fig:h2bounds} shows various bounds and values for the
  monomer-dimer entropy $h_2(p)$ of dimer density $p \in [0,1]$ in
  the $4$-regular $2$-dimensional grid. FT is the Friedland-Tverberg
  lower bound $fh_4(p)$ of (\ref{FTLowerBound}), h2 is the true
  monomer-dimer entropy equal to $\max_{p \in [0,1]}h_2(p)$ (it is
  known to a precision much greater than the picture resolution). The
  crosses marked B are Baxter's computed values \cite{Bax1}. ALMC is
  the function $gh_4(p)$ of (\ref{defghrp}), conjectured to be a
  lower bound in the Asymptotic Lower Matching Conjecture. AUMC is
  the monomer-dimer entropy $h_K(p)$ of dimer density $p$ in a
  countable union of $K_{4,4}$, given by (\ref{presK})--(\ref{pt})
  and conjectured to be an upper bound by the Asymptotic Upper
  Matching Conjecture. Notice that AUMC goes a little over h2: a
  countable union of $K_{4,4}$ has a higher monomer-dimer entropy
  than an infinite planar grid.

\begin{figure}[th]
    \begin{center}
    \includegraphics[width=0.75\textwidth]{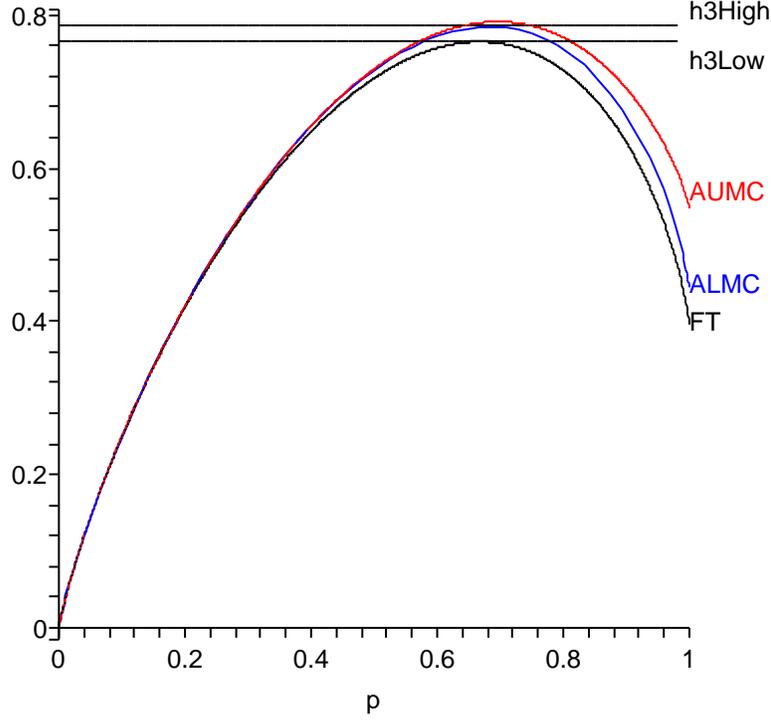}
    \caption{Monomer-dimer tiling of the $3$-dimensional grid:
    entropy as a function of dimer density. FT is the
    Friedland-Tverberg lower bound, h3Low and h3High are the known
    bounds for the monomer-dimer entropy. ALMC is the Asymptotic
    Lower Matching Conjecture. AUMC is the entropy of a countable
    union of $K_{6,6}$, conjectured to be an upper bound by the
    Asymptotic Upper Matching Conjecture.}\label{Fig:h3bounds}
    \end{center}
  \end{figure}

 Figure~\ref{Fig:h3bounds} shows similarly various bounds and
  values for the monomer-dimer entropy $h_3(p)$ of dimer density $p
  \in [0,1]$ in the $6$-regular $3$-dimensional grid. FT is the
  Friedland-Tverberg lower bound $fh_6(p)$ of (\ref{FTLowerBound}),
  h3Low and h3High are the best known lower and upper bounds for the
  true monomer-dimer entropy equal to $\max_{p \in [0,1]}h_3(p)$.
  ALMC is the function $gh_6(p)$ of (\ref{defghrp}), conjectured to
  be a lower bound in the Asymptotic Lower Matching Conjecture. AUMC
  is the monomer-dimer entropy $h_K(p)$ of dimer density $p$ in a
  countable union of $K_{6,6}$, given by (\ref{presK})--(\ref{pt})
  and conjectured to be an upper bound by the Asymptotic Upper
  Matching Conjecture. Notice that AUMC goes a little over h3High: a
  countable union of $K_{6,6}$ has a higher monomer-dimer entropy
  than an infinite cubic grid.

 \section{A proof of some case of the ALMC}
 In this section we prove the following case of ALMC:
 \begin{theo}\label{palmcc}  Let $r\ge 3$ be an integer.  Then
 the asymptotic lower matching conjecture (\ref{lasmbdfp})
 holds for $p_s=\frac{r}{r+s}, s=0,1,2,\ldots$.
 \end{theo}

 The proof of this theorem follows from the following results.
 \begin{theo}\label{frtverv1}  Let $p:R^n \to \R$ be  positive
 hyperbolic of degree $m\in [1,n)$.
 Assume that $\deg_i p\le r_i \in [1,m]$ for $i=1,\ldots, n$.
 Rearrange the sequence $r_1,\ldots,r_n$ in an increasing
 order $1\le r_1^*\le r_2^*\le \ldots \le r_n^*$.  Let $s\in\N$.
 Let $k\in [1,n]$ be
 the smallest integer such that $r_k^*+s> n-k$.
 Then
 \begin{eqnarray}
 \sum_{1\le i_1 <\ldots<i_m\le n} \frac{\partial ^m p}{\partial
 x_{i_1}\ldots\partial x_{i_m}} (\0) \ge \nonumber\\
 \frac{(s n)!}{s^{n-m}(n-m)!((s-1)n+m)!}\frac{(n-k+1)!}
 {(n-k+1)^{n-k+1}}\prod_{j=1}^{k-1}
 (\frac{r_j^*+s-1}{r_j^*+s})^{r_j^*+s-1} \Cp p.
 \label{frtver1v1}
 \end{eqnarray}
 \end{theo}

 \proof  Let $q:\R^n \to \R$ be positive
 hyperbolic of degree $n-m$ with $\deg_i q\le s$ for $i=1,\ldots,n$
 and $\Cp q=1$.  Then $f=pq:\R^n \to \R$ is positive hyperbolic of degree $n$ with
 $\Cp f\ge\Cp p$ and $\deg_i f\le r_i+s$ for $i=1,\ldots,n$.
 Apply Theorem \ref{frtver} to $f$ to deduce
 \begin{eqnarray}
 \sum_{1\le i_1<\ldots<i_m\le n}
 \frac{\partial ^m p}{\partial
 x_{i_1}\ldots\partial x_{i_m}} (\0)\;
 \frac{\partial ^{n-m} q}{\partial
 x_{i'_1}\ldots\partial x_{i'_{n-m}}} (\0) \ge \nonumber\\
 \frac{(n-k+1)!}{(n-k+1)^{n-k+1}}\prod_{j=1}^{k-1}
 (\frac{r_j^*+s-1}{r_j^*+s})^{r_j^*+s-1} \Cp p, \label{frtver2v1}
 \end{eqnarray}
 where $1\le i'_1<\ldots<i_{n-m}'\le n$  and
 $\{i_1,\ldots,i_m,i_1',\ldots,i_{n-m}'\}=\an{n}$.

 Let $A\in \Delta(n,s)$ and choose
 $q={n\choose n-m}^{-1}p_{n-m,\frac{1}{s}A}(\x)$ as in
 (\ref{exhpk}).    Note
 $$\frac{\partial ^{n-m} q}{\partial
 x_{i'_1}\ldots\partial x_{i'_{n-m}}} (\0)=
 \frac{1}{{n\choose n-m}s^{n-m}}\perm
 _{n-m} A[\an{n}|J'].$$
 Now take the expected value of the left-hand side of the
 inequalities (\ref{frtver2v1}) corresponding to all
 $A\in\Delta(n,s)$.
 Use Corollary \ref{mexpcvalcor} to deduce that the coefficient of
 each $\frac{\partial ^m p}{\partial
 x_{i_1}\ldots\partial x_{i_m}} (\0)$ is $\frac{s^{n-m}(n-m)!((s-1)n+m)!}{(s
 n)!}$.  \qed

 \begin{corol}\label{frtverv2}  Let $p:R^n \to \R$ be  positive
 hyperbolic of degree $m\in [1,n)$.
 Assume that $\deg_i p\le r \in [1,m]$ for $i=1,\ldots, n$.
 Let $s\in\N$ and $k=n-r-s+1\ge 1$.
 Then
 \begin{eqnarray}
 \sum_{1\le i_1 <\ldots<i_m\le n} \frac{\partial ^m p}{\partial
 x_{i_1}\ldots\partial x_{i_m}} (\0) \ge \nonumber\\
 \frac{(s n)!}{s^{n-m}(n-m)!((s-1)n+m)!}\frac{(r+s)!}{(r+s)^{r+s}}
 \big(\frac{r+s-1}{r+s})^{(r+s-1)(n-r-s)} \Cp p.
 \label{frtver1v2}
 \end{eqnarray}
 \end{corol}

 \begin{theo}\label{frtverv3}  Let $B\in\R_+^{n\times n}$ be a
 doubly stochastic matrix with at most $r$ nonzero entries in each
 column.
 Let $s\in\N$ and $k=n-r-s+1\ge 1$.
 Then for each $m\in [1,n)$
 \begin{equation}
 \perm_m B\ge
 \frac{(s n)!{n\choose m}}{s^{n-m}(n-m)!((s-1)n+m)!}\frac{(r+s)!}{(r+s)^{r+s}}
 \big(\frac{r+s-1}{r+s})^{(r+s-1)(n-r-s)}.
 \label{frtver1v3}
 \end{equation}

 \end{theo}

 \proof  Let $p=p_{m,B}(\x)$ as defined by (\ref{exhpk}).
 Then (\ref{frtver1v3}) follow from Corollary
 \ref{frtverv2}.  \qed

 Let $G_n\in\Gamma(n,r)$.  Then $G_n$ is represented by its incidence matrix
 $A_n\in \Delta(n,r)$.  Let $B_n:=\frac{1}{r} A_n$.  Then $B_n$ is a
 doubly stochastic matrix where each row and column of $B_n$ has at
 most $r$ positive entries.  Clearly, the ALMC conjecture follows from
 the following stronger conjecture:

 \begin{con}[The Asymptotic Lower $r$-Permanent Conjecture ]
      \label{lasrpc} \mbox{}\newline
      For $r \geq 2$, let $B_n, n=1,2,\ldots$ be a
      sequence of $n\times n$ doubly stochastic matrices, where each
      column of each $B_n$ has at most $r$-nonzero entries.
      Let $k_n \in [0,n]$, $n = 1,2,\ldots$
      be a sequence of integers with $lim_{n \to \infty}
      \frac{k_n}{n} = p\in (0,1]$.  Then
      \begin{equation}\label{lasrppcp}
       \liminf_{n \to \infty} \frac{\log \perm_{k_n} B_n}{2n} \geq gh_r(p) -\frac{p}{2}\log r.
      \end{equation}
 \end{con}

 Theorem  \ref{palmcc} follows from the following result:

 \begin{theo}\label{palrpc}  Let $r\ge 3, s\ge 1$ be integers.
 Let $B_n, n=1,2,\ldots$ be a
      sequence of $n\times n$ doubly stochastic matrices, where each
      column of each $B_n$ has at most $r$-nonzero entries.
      Let $k_n \in [0,n]$, $n = 1,2,\ldots$
      be a sequence of integers with $lim_{n \to \infty}
      \frac{k_n}{n} = p\in (0,1]$.  Then
      \begin{eqnarray}\label{lasrppbp}
      \liminf_{n \to \infty} \frac{\log \perm_{k_n} B_n}{2n} \geq
      \frac{1}{2}\left(-p\log p - 2(1-p)\log (1-p)\right) +\\
       \frac{1}{2}\left(
      (r+s-1)\log(1-\frac{1}{r+s})-(s-1+p)\log(1-\frac{1-p}{s})
       \right).\nonumber
      \end{eqnarray}
 Moreover,
 the Asymptotic Lower $r$-Permanent Conjecture \ref{lasrpc}
 holds for $p_s=\frac{r}{r+s}, s=0,1,2,\ldots$.
 \end{theo}

 \textbf{Proof of Theorem \ref{palrpc}}.
 Apply the inequality (\ref{frtver1v3}) to $B_n$ for $m=k_n$.  Take
 the logarithm of the both sides of this inequality and let
 $n\to\infty$.  A straightforward calculation for the right-hand
 side, using the Stirling's formula,  yields the inequality (\ref{lasrppbp}).
 Assume that $p=p_s=\frac{r}{r+s}$.  Then
 $\frac{1-p_s}{s}=\frac{1}{r+s}=\frac{p_s}{r}$.  Then the
 right-and side of (\ref{lasrppbp}) is equal to
 $gh_r(p_s) -\frac{p_s}{2}\log r$.  Hence
 the asymptotic lower $r$-permanent conjecture \ref{lasrpc}
 holds for $p_s=\frac{r}{r+s}, s=1,2,\ldots$.

 We now discuss the case $s=0$, i.e.
 $p=p_0=1$.
 Let $B=(b_{ij})_{i,j=1}^n$ be any $n\times n$ nonnegative matrix.  Denote
 by $G(B)=(V,E)$ the bipartite graph induced by $B$, i.e. the edge
 $(i,j)$ is in $E$, if and only if $b_{ij}>0$.
 Then $B$ induces the weighted graph on $G$, where the weight
 of the edge $(i,j)$ is $b_{ij}$.
 Let
 $p_B(x)=x^n+\sum_{m=1}^n (-1)^m \perm_m (B)$.
 $p_B(x)$ is called the matching polynomial of the weighted graph $G$.
 Heilmann and Lieb showed in \cite{HL} that $p_B(x)$
 has nonnegative roots.  (See also \cite{LP}.)
 Hence the arithmetic-geometric inequality for the
 elementary symmetric polynomials of the nonnegative roots of $p_B(x)$
 yields the inequality
 $\perm_{m} B\ge {n \choose m} (\perm B)^{\frac{m}{n}}$.
 (See \cite{Wan} for the case of $m$-matchings in bipartite graphs.)

 Use Theorem \ref{gurschr} to deduce that
 $\perm B_n\ge \frac{r!}{r^r}\big(\frac{r}{r-1}
 \big)^{r(r-1)}\big(\frac{r-1}{r}\big)^{(r-1)n}$.
 Apply the above two inequalities for the sequence $B_n$ and $m=k_n$
 for $n=1,2,\ldots$
 to deduce the case $p=1$.   \qed

 Let $d=6$ and
 $p_3:=\frac{6}{9}=\frac{2}{3}$.  Then Theorem \ref{palmcc} yields
 that $h_3(\frac{2}{3})\ge .7845241927$, which implies that
 $h_3=\max_{p\in[0,1]} h_3(p)\ge  .7845241927$.  This improves the lower
 bound implied by (\ref{lasmbdfp}) $h_3\ge .7652789557$
 \cite{FP}.  The computations in \cite{FP} yield that $h_3\le
 .7862023450$.  Thus $h_3\in [.7845, .7863]$.

 \section{Matching in general graphs - hafnians}

 Let $G=(V,E)$ be a graph on the set of vertices $V$ and the set of
 edges $E$.  Assume that $\#V=N$.  Then $G$ is represented by a
 symmetric $0-1$ matrix $B=B(G)$ with $0$ diagonal.  If $G$ has a perfect
 matching then $N=2n$ is even.  If $G$ is bipartite and $V=V_1\cup V_2$,
 where $V_1=\{1,\ldots,n\}, V_2=\{n+1,\ldots,2n\}$, we deduce that
 \begin{equation}\label{biparblmat}
 B=\left(\begin{array}{cc}\0&A\\A\trans&\0\end{array}\right),
 \end{equation}
 where
 $A$ is the representation matrix of the bipartite graph $G$.
 As explained above, the number of $m$-matching in the bipartite
 graph $G$ is $\perm _m A$.

 In this section we discuss the $m$-matching of a general graph $G$,
 and the related function $\haff_m B$ which counts the number of
 $m$-matching in $G$.
 Let $\alpha,\beta\subseteq \{1,2,\ldots,N\}$ be two nonempty sets
 of cardinality $i,j$, $\#\alpha=i,\#\beta=j$, respectively.
 We then arrange the elements of
 $\alpha=\{\alpha_1,\ldots,\alpha_i\}$ and
 $\beta=\{\beta_1,\ldots,\beta_j\}$ in an increasing order:
 $1\le\alpha_1<\ldots<\alpha_i\le N, 1\le\beta_1<\ldots<\beta_j\le
 N$.  For $B=(b_{st})\in \C^{N\times N}$ we denote by $B[\alpha|\beta]\in
 C^{i\times j}$ the submatrix $(b_{\alpha_s\beta_t})_{s,t=1}^{i,j}$.
 Denote by $\rS_l(\R)\supset \rS_l(\R_+), \rS_{l,0}(\R)\supset\rS_{l,0}(\R_+)$
 the space of real valued $l\times l$ symmetric matrices, the cone
 $l\times l$ symmetric matrices with nonnegative entries, the
 subspace of real valued $l\times l$ symmetric matrices with zero
 diagonal, and the subcone of $l\times l$ symmetric matrices with zero
 diagonal and nonnegative entries respectively.
 Let $B\in \rS_{N}(\R)$ and an integer $m\in [1,\frac{N}{2}]$.  Then the $m
 -th$ \emph{hafnian} of $B$ is defined as
 \begin{equation}\label{defhafm}
 \haff _m B=2^{-m}\sum_{\alpha,\beta\subset\{1,\ldots,N\},
 \#\alpha=\#\beta=m, \alpha\cap\beta=\emptyset} \per B[\alpha,\beta].
 \end{equation}
 That is if $(i_1,j_1),\ldots, (i_m,j_m)$ is an $m$ matching of a
 complete graph $K_{N}$ on $N$ vertices, then the product $b_{i_1j_1}\ldots b_{i_mj_m}$
 appears exactly once in $\haff_m B$.  Since $b_{i_lj_l}=b_{j_li_l}$
 there are $2^m$ choices of $\alpha$ and $\beta$ for which this
 product appear, we need to use the factor $2^{-m}$ in the above
 definition of $\haff_m B$.  If $B=B(G)$ then $\haff_m B$ gives the
 number of $m$-matching in $G$.  Note that from the definition
 of $\perm _m B$ it follows that $\haff_m B\le 2^{-m}\perm_m B$.
 Equivalently, it is straightforward to show:
 \begin{equation}\label{hafparder1}
 \haff_m B=(2^{m} m!)^{-1} \sum_{1\le i_1 <\ldots<i_{2m}\le N}
 \frac{\partial ^{2m} }{\partial
 x_{i_1}\ldots\partial x_{i_{N}}} (\x\trans B\x)^m, \quad B\in \rS_{N}(\R) .
 \end{equation}
 Unfortunately, the quadratic polynomial $\x\trans B\x$ is not always positive hyperbolic.
 Note that $\haff_m B$ does not depend on the value of the diagonal
 entries $B$.  Let $B^{(0)}$ be the matrix obtained from by
 replacing the diagonal entries of $B$ by zero elements.
 Then $\haff_m B=\haff_m B^{(0)}$.

 For $B\in \rS_l(\R)$ we denote by $\lambda_1(B)\ge\ldots\ge
 \lambda_l(B)$ the $l$ eigenvalues of $B$ counted with their
 multiplicities and arranged in the decreasing order.
 As usual for $B,C\in \rS_l(\R)$ we denote $C\succeq B$ if $C-B$ is a
 nonnegative definite matrix.  The maxmin, or minmax charaterization
 of $\lambda_k(B)$ yields that if $C\succeq B$ then $\lambda_k(C)\ge
 \lambda_k(B)$ for $k=1,\ldots,l$.  In particular, if
 $B\in\rS_l(\R)$ has nonnegative diagonal entries then $B\succeq B^{(0)}$
 and $\lambda_k(B)\ge \lambda_k(B^{(0)})$ for $k=1,\ldots,l$.

 The following result is well known and we bring
 its proof for completeness.
 \begin{lemma}\label{quadfrmph}  Let $B\in \rS_n(\R)$ and $n\ge 2$.  Then
 $\x\trans B\x$ is positive hyperbolic if and only if
 $\0\ne B\in\rS_n(\R_+)$
 and $0\ge \lambda_2(B)$.
 \end{lemma}
 \proof  Assume that $\x\trans B\x$ is positive hyperbolic.
 Then Proposition \ref{nonnegcphp} yields that
 $B\in\rS_n(\R_+)$.  Since $\x\trans B\x>0$ for
 each $\x>\0$ $B\ne \0$.  Hence $\lambda_1(B)>0$.

 Observe next that the positive
 hyperbolicity  of $\x\trans B\x$ is equivalent to
 \begin{equation}\label{quadfrmph1}
 (\x\trans B\y)^2\ge (\x\trans B\x)(\y\trans B\y), \textrm{ for
 any }
 \x>\0, \y\in\R^n.
 \end{equation}
 Clearly, the above condition holds for any $\x\ge \0$.
 The Perron-Frobenius theorem yields that
 there exists $\0\ne \x\in \R^n_+$ such that
 $B\x=\lambda_1(B)\x$.  Let $\0\ne \y\in \R^n$ and assume that
 $\y\trans\x=0$.  Then (\ref{quadfrmph1}) yields that $\y\trans B\y\le 0$.
 Hence $\lambda_2(B)\le 0$.

 Vice versa suppose that $\0\ne B\in \rS_n(\R_+)$
 and $\lambda_2(B)\le 0$.
 Recall that there exists
 a permutation matrix $P\in \{0,1\}^{n\times
 n}$ such that $P\trans BP$ is a block diagonal matrix
 $\diag(B_1,\ldots,B_k)$, where $B_i\in
 \rS_{n_i}( \R_+)$ is irreducible.  So $\lambda_1(B_i)>0$ unless
 $n_i=1$ and $B_i=0$.  Hence our assumptions yield that we may
 assume that $\0\ne B_1\in \rS_{n_1}(\R_+)$ is irreducible and
 $B_i=0$ for $i=2,\ldots,k$.  Thus it is enough to show that
 $\x\trans B\x$ is positive hyperbolic for an irreducible $B\in
 \rS_n(\R_+)$, where $n\ge 2$.  Clearly $\x\trans B\x>0$ for $\x>0$.
 Thus it is left to show that (\ref{quadfrmph1}) holds.
 Assume first that
 $B\x=\lambda_1(B)\x, \x >0$.  Then (\ref{quadfrmph1}) follows straightforward.
 Suppose $\x>0$ is any vector.  Then there exists a unique diagonal
 matrix $D$, with positive diagonal entries such that
 $B\x=D^{-2}\x$.  That is $DBD(D^{-1}\x)=D^{-1}\x$.  Replacing
 $B,\x,\y$ by $DBD, D^{-1}\x,D^{-1}\y$ we deduce the inequality
 second inequality of (\ref{quadfrmph1}).  \qed

 \begin{defn}\label{kpartdef}  Let $G=(V,E)$ be a graph on the set
 of vertices $V$.  Then for $k\ge 2$ $G$ is called $k$-partite,
 if the following condition holds.  There exists a decomposition
 of $V$ to a disjoint union of $k$ nonempty sets $V_1,\ldots,V_k$
 such that each edge $e\in E$ connects $V_i$ to $V_j$ for some $i\ne
 j$.
 $G$ is called a complete $k$-bipartite, if there exists a decomposition
 of $V$ to a disjoint union of $k$ nonempty sets $V_1,\ldots,V_k$
 such that $E$ consists of all edges from  $V_i$ to $V_j$ for all $1\le
 i < j \le k$.

 \end{defn}

 Note that a complete graph $G$ with $n\ge 2$ vertices is complete $n$-
 partite.\\
 We need the following elementary result whose proof is straightforward.
 \begin{prop}\label{3-minor}
 Let $n \geq 3$ and $F=[f_{ij}]\in \rS_n(\R)$
 with $1's$ on the main diagonal.
 Suppose that
 and all subsets of cardinality three $\alpha = \{1 \leq i< j <k \leq n\}$
 the principal submatrix $F[\alpha,\alpha]$ is
 nonnegative definite . If $f_{ij} = 1 , f_{jk} =1$
 then $f_{ik}=1$.
 \end{prop}

 \begin{lemma}\label{01pathypq}  Let $B=[b_{ij}]\in\rS_{n}(\R_+)$, i.e.
 $B$ is a real $n\times n$ symmetric matrix with nonnegative entries.
 Denote by $G(B)=(\lan n\ran, E)$ the graph, (with no self-loops), induced by
 $B$, i.e. $(i,j)\in E$ if and only if $i\ne j$ and $b_{ij}>0$.
 Assume that $B$ is
 irreducible and $\lambda_2(B)\le 0$.  Then $G(B)$ is a complete
 $k$-partite graph for some $k\in [2,n]$.
 \end{lemma}

 \proof
 Since $B^{(0)}$ is irreducible and $B\succeq B^{(0)}$, it is enough
 to prove the lemma in the case where all the diagonal entries of
 $B$ are equal to zero, i.e. $B=B^{(0)}$.
 Let $\x=(x_1,\ldots,x_n)\trans$ be the unique positive eigenvector of
 $B$, corresponding to the maximal eigenvalue $\lambda_1(B)$,
 satisfying the condition $\x\trans \x=\lambda_1(B)$.
 Then $A = \x \x^{T} - B$
 and $A\succeq\0 $. Let $D = \diag(x_{1},...,x_{n})$.
 Then the zero pattern of the matrix $C=D^{-1}BD^{-1}$ is
 the same as of the matrix $B$.  Let
 $F=[f_{ij}]: = \1 \1\trans - C$ , where $\1:=(1,\ldots,1)\trans$
 is the vector of all ones.  Then $F=[f_{ij}] \succeq \0$.
 Notice that $f_{ii} = 1 ,1 \leq i \leq n$, and for $i\ne j$
 $f_{ij}=1$ if and only $b_{ij} = 0$, (i.e. the vertices $i,j$
 are not connected in the graph $G(B)$).
 As any principal submatrix of $F$ is nonnegative definite,
 it follows from Proposition \ref{3-minor} that for any triplet
 $1 \leq i < j < k \leq n$
 such that if $f_{ij} = 1 , f_{jk} =1$ then also $f_{ik}=1$.
 In other words, the relation
 $"i \thicksim j" \Longleftrightarrow "b_{ij} = 0" $
 is an equivalence relation.
 Therefore the graph $G(B)$ is complete $k$-partite, where each
 equivalence class of vertices corresponds to some class $V_i$
 in the $k$-partite graph.

 \begin{theo}\label{hypkpgrh}  Let $G$ be a graph on $n >1$ vertices,
 and denote by $A(G)$ the incidence matrix of $G$.
 Then $\x\trans A(G)\x$ is positive hyperbolic if and only if
 $G$ is a union of a complete $k(\ge 2)$-partite graph on at least two
 vertices and of isolated vertices.
 \end{theo}
 \proof  We first show that for a complete $k$-partite graph $G$ on at least
 $n\ge 2$ vertices $\lambda_2(G)\le 0$, which is equivalent to the positive
 hyperbolicity of $\x\trans A(G)\x$ in view of Lemma \ref{quadfrmph}.
 Let $J_n$ be a symmetric matrix whose all entries are equal to $1$.
 Then $J$ is rank one matrix with $\lambda_1(J_n)=n$ and
 $\lambda_i(J_n)=0$ for $i=2,\ldots,n$.
 Let $k\in [1,n]$, $n_i\in\N, i=1,\ldots,k, 1\le n_k\le \ldots\le n_1$
 and $n_1+\ldots+n_k=n$.
 Consider the block diagonal matrix
 $J(n_1,\ldots,n_k):=\diag(J_{n_1},\ldots,J_{n_k})$.
 Clearly $J(n_1,\ldots,n_k)$ is a nonnegative definite matrix.
 It is
 straightforward to see that renaming the vertices of $G$, we will
 obtain that $A(G)=J_n-J(n_1,\ldots,n_k)$ for some unique
 $n_1\ge\ldots n_k\ge 1$.  Then minimax characterization of
 $\lambda_2(A(G))$ yields that $\lambda_2(A(G))\le \lambda_2(J_n)=0$.
 Hence $\x\trans A(G)\x$ is positive hyperbolic.

 Assume now that $\x\trans A(G)\x$ is positive hyperbolic.
 Therefore $G$ must have at least
 one edge and $\lambda_2(A(G))\le 0$.
 Hence $G$ has one connected component
 containing at least two vertices and a union of
 isolated vertices.
 Without loss of generality we assume that $G$ is
 connected.  Then $A(G)$ satisfies the assumptions
 of Lemma \ref{01pathypq} and $G$ is $k$-partite.
 \qed

 \begin{rem}\label{reverse}
 Let $A$ be a symmetric $n \times n$ matrix with nonnegative entries.
 It is straightforward to show that the polynomial $\x\trans A \x$ is positive hyperbolic
 if and only if the function $\sqrt{\x\trans A \x}$ is concave on the positive
 orthant $R^{n}_{+}$.
 If $A$ is real, symmetric, nonnegative definite then $\sqrt{\x\trans A \x}$
 is convex on $R^n$.
 In view of Remark \ref{concave} , it is natural to conjecture that
 if $A$ is
 $2n \times 2n$  real, symmetric, nonnegative definite then the reverse van der Waerden bound holds :
 $$
 2^n n! \haff_n B \leq \frac{(2n)!}{(2n)^{2n}} \Cp (p) ,\quad p(\x): = (\x\trans A \x)^n .
 $$

 \end{rem}

 \section{Algorithmic applications}

 One of the main purposes of this paper is
 to construct a generating homogeneous polynomial
 $p(\x)$ of degree $n$,
 for some quantity of interest
 $Q$ as hafnian: $\haff A$, sums of subhafnians: $\haff_{m} A$,
 permanent: $\per A$, sum of
 subpermanents: $\per_m A$,  such that $Q= \frac{\partial
 ^n p}{\partial x_1 ....\partial x_n} (\0)$. If such polynomial
 is positive hyperbolic then we can apply the results
 from \cite{Gur1} and the results of the previous sections
 of this paper to get a lower bound on $Q,$ and even to get
 deterministic polynomial-time algorithms to approximate $Q$ within
 simply exponential factor as in \cite{Gur1}. In the general, (not positive
 hyperbolic case), we can use this representation to obtain exact
 algorithm to compute $Q$ in $2^{n} poly(n)$ number of arithmetic
 operations provided the the generating polynomial $p$ can be itself
 evaluated in $poly(n)$ number of arithmetic operations. We present
 below some examples of this approach.

 \subsection{Formula for $\sum_{1\le i_1 <\ldots<i_m\le n} \frac{\partial ^m p}{\partial
 x_{i_1}\ldots\partial x_{i_m}} (\0)$ }

 Our exact algorithms are based on the following
 elementary identity \ref{ryzer}.

 Let $p(\x)$ be a polynomial of degree $m$ in $n \geq m$ variables , $p(0) =0$.
 Define
 \begin{equation}\label{defS_{i}}
 s_{i} = \sum_{b_{j} \in \{0,1\} , \sum{1 \leq j \leq m} = i} p(b_1,...,b_m) .
 \end{equation}

 Let $\bd_{n} = (d_{n,1},...,d_{n,m-1})$ be the unique
 solution of the system of linear equations $\bd_n A = (-1,...,-1)$ ,
 where the $m-1 \times m-1$ lower triangular matrix $A=[a_{ij}]$ is defined as follows:
 $$a_{ij} = {n-j \choose i-j} \textrm{ if }  i \geq j \textrm{ and }
 a_{ij}=0 \textrm{ otherwise}.$$
 Then the following equality holds :
 \begin{equation}\label{ryzer}
 \sum_{1\le i_1 <\ldots<i_m\le n} \frac{\partial ^m p}{\partial
 x_{i_1}\ldots\partial x_{i_m}} (\0) = p(1,...,1) + \sum_{1\le j \leq m-1} s_{j}d_{n,j}
 \end{equation}
 Notice that this formula requires $\sum_{0 \leq j \leq n-1} {n \choose j}$
 evaluations of the polynomial $p$.

 The formula \ref{ryzer} follows from the following obvious identities :
 $$
 s_{i} = \sum_{1 \leq j \leq i} {n-j \choose i-j} c_{j} , 1 \leq i \leq m ,
 $$
 where $c_{j}$ is the sum of the coefficients of all monomials in $p$ involving
 exactly $j$ variables.

 The formula \ref{ryzer} is, in a sense, optimal for $n=m$, i.e., if for all homogeneous
 polynomials $p(\x)$ of degree $n$
 $$
 \frac{\partial^n p}{\partial x_1... \partial x_n} (\x) =
 \sum_{1 \leq i \leq k} a_i p(\z_i), \quad a_i \in \C , \z_i \in \C^n ,
 $$
 then $k \geq {n \choose \frac{n}{2} } \approx
 \frac{2^{n}}{\sqrt{n}}$ \cite{Gur-Hyp}.

 \subsection{Ryser' like formulas for sums of subhafnians and subpermanents}
 \begin{enumerate}


 \item
 Let $B\in \rS_{N}(\R),\x:=(x_1,\ldots,x_{N})\trans \in \C^{N}$ and $m\in
 [ 1,n]\cap\N$.
 Define, as in the proof of Theorem 3.1, the polynomial :
 $$
 P_{m}(\x) = \frac{1}{ 2^{m}m!}(\x\trans B\x)^m
 $$

 As  $\haff_m B=(2^{m} m!)^{-1} \sum_{1\le i_1 <\ldots<i_{2m}\le N}
 \frac{\partial ^{2m} }{\partial
 x_{i_1}\ldots\partial x_{i_{2m}}} (\x\trans B\x)^m$ , the application of
 \ref{ryzer} gives
 the following Ryser-like formula for $\haff_m B$:
 \begin{equation}\label{ryzer-haf}
 \haff_m B= P_{m}(1,...,1) + \sum_{1\le j \leq m-1} s_{j}d_{n,j},
 \end{equation}
 where $s_i$ are defined by (\ref{defS_{i}}) for $p=P_m$.
 The formula (\ref{ryzer-haf}) provides
 $SB(N,m) (O(N^{2}) + O(\log (m)))$ , algorithm to compute $\haff_{m}
 B$, where
 $SB(N,m): =\sum_{0 \leq j \leq 2m-1} {N \choose j}$ .
 \item
 For $\x:=(x_1,...,x_n)\trans
 \in \C^n$ let
 $$
 S_m(\x) = \sum_{1 \leq i_{1} <...<i_{m} \leq n}
 x_{i_{1}}...x_{i_{m}}.
 $$
 be the $m-th$ symmetric function of $\x$.
 Let $A$ be $n \times n$ complex matrix and define
 $p_{m}(\x):= S_{m}(A \x)$.  Then
 $p_{m}(\x)$ can be evaluated in $O(n^2)$ arithmetic
 operations and
 $$
\per_{m} A = \sum_{1\le i_1 <\ldots<i_{m}\le n}
 \frac{\partial ^{m} }{\partial
 x_{i_1}\ldots\partial x_{i_{m}}} p_{m}(\x).
 $$
 Which gives the following algorithm to evaluate $ \per_{m} A$
 \begin{equation}\label{ryzer-perm}
 \per_{m} A =p_{m}(1,...,1) + \sum_{1\le j \leq m-1} s_{j}d_{n,j},
 \end{equation}
 where $s_i$ are defined by (\ref{defS_{i}}) for $p=p_m$.
 The formula \ref{ryzer-perm} provides $SB(N,m) (O(N^{2}))$ algorithm to compute $\per_{m} A$.\\

 Notice that the naive algorithm , i.e. computing and adding all $m \times m$ subpermanents,
 requires $( {n \choose m})^2 2^{m} O(m)$ arithmetic operations .

 \end{enumerate}

 \subsection{Positive hyperbolic polynomials and convex relaxations}

 In this section we always assume that
 $\x=(x_1,\ldots,x_n)\trans, \y=(y_1,\ldots,y_n)\trans\in\C^n$,
 $\1=(1,\ldots,1)\trans$.
 Suppose that a positive hyperbolic polynomial
 $$
 p(\x) =\sum_{ \sum_{1 \leq i \leq n} r_{i} =n }
 a_{(r_1,...,r_n) } \prod_{1 \leq i \leq m} x_{i}^{r_{i}} , m \geq n
 $$
 has nonnegative integer coefficients and is given as an oracle.
 I.e. we don't have a list coefficients , but can evaluate
 $p(\x)$ on rational inputs . The number $\log p(\1)$
 measures the complexity of the polynomial $p$ .

 A deterministic polynomial-time oracle algorithm is any algorithm
 which evaluates the given polynomial $p(\x)$ at a number of rational
 vectors $\q^{(i)}=(q_{1}^{(i)},...,q_{n}^{(i)})$ which is polynomial
 in $n$ and $\log p(\1)$; these rational vectors $\q^{(i)}$ are
 required to have bit-wise complexity which is polynomial in $n$ and
 $\log p(\1)$; and the number of additional auxiliary
 arithmetic
 computations is also polynomial in $n$ and $\log p(\1)$.\\
 If the number of oracle calls, (evaluations of the given polynomial
 $p(\x)$), the number of additional auxiliary arithmetic computations
 and bit-wise complexity of the rational input vectors $\q^{(i)}$
 are all polynomial in $n$, (no dependence on $\log p(\1)$)
 then such algorithm is called  deterministic strongly polynomial-time
 oracle algorithm.\\
 The following theorem was proved in \cite{Gur1}.

 \begin{theo}\label{Conv}
 There exists a deterministic polynomial-time oracle algorithm, which
 computes for given as an oracle indecomposable positive hyperbolic
 polynomial $p(\x)$ a number $F(p)$, satisfying the
 inequality
 $$
 \frac{\partial^n p}{\partial x_1...\partial x_n} (\0) \leq F(p)
 \leq 2\frac{n^n}{n!}\frac{\partial^n p}{\partial x_1...\partial x_n}
 (\0)
 $$

 \end{theo}

 Our goal in this paper is to extend Theorem 7.1 to approximate
 \par\noindent
 $\sum_{1\le i_1 <\ldots<i_m \le n} \frac{\partial ^m q}{\partial
  x_{i_1}\ldots\partial x_{i_m}} (\0)$ for a given positive
  hyperbolic polynomial $q(\x)$ of degree $m$ .

  The algorithm behind Theorem 7.1 is based on two observations . First,
  \begin{equation}\label{conv-rel}
  \log \Cp p = \inf_{\sum_{1 \leq i \leq n} y_I = 0} \log
  p(e^{\y}),
 \end{equation}
 where we used the notation of Proposition \ref{logconsumexp}.
 If the coefficients of the polynomial $p$ are nonnegative then the
 functional $\log p(e^{\y})$ is convex ,
 the indecomposabilty of the polynomial $p$ is exactly uniqueness
 and existence of the minimum in \ref{conv-rel} . \\
 Second point is the inequality :
 \begin{equation}\label{cap-waer}
 \Cp p\frac{n!} {n^n}   \leq    \frac{\partial^n p}{\partial
 x_1...\partial x_n} (\0) \leq \Cp p
 \end{equation}
 If the positive hyperbolic polynomial $p$ is not indecomposable , we
 need first to check if $\Cp p > 0$ . If $\Cp p = 0$ then also
 $\frac{\partial^n p}{\partial x_1...\partial x_n} (\0) =0$ .
 In the case $\Cp p > 0$ we "slightly" perturb the polynomial $p$
 to get the indecomposability.\\
 Theorem 3.1 in this paper provides an analogue of the left
 inequality in (\ref{cap-waer}) for positive hyperbolic polynomial
 $q(\x)$ of degree $m < n$ . The problem is that in this
 case we don't have the right inequality. I.e. it is possible that
 $\sum_{1\le i_1 <\ldots<i_m \le n} \frac{\partial ^m}{\partial
  x_{i_1}\ldots\partial x_{i_m}} q(\0) > 0$ but $\Cp(q) = 0$ .
  This problem can be easily overcome
  by the following equivalent reformulation of Theorem 3.1:

 \begin{theo}\label{Tver-new}
 Consider a positive hyperbolic polynomial $q(\x)$ of degree
 $m < n$ . Define a positive hyperbolic polynomial $p(\x) =
 q(\x)(\frac{\sum_{1 \leq i \leq n} x_{i}}{n})^{n-m}$ and
 $$
 D_{m}(q) = \frac{ (n-m)!}{n^{n-m}} \sum_{1\le i_1 <\ldots<i_m \le n}
 \frac{\partial ^m q}{\partial
  x_{i_1}\ldots\partial x_{i_m}} (\0).
 $$
 Then the following inequality holds :
 \begin{equation}\label{waer-new}
 \frac{n!} {n^n} \Cp p \leq D_{m}(q) = \frac{\partial^n p}{\partial
 x_1...\partial x_n} (\0)  \leq \Cp p .
 \end{equation}
 \end{theo}

 \begin{corol}\label{zer--check}
 Let $q( \x)$ be a positive hyperbolic polynomial of degree
 $m < n$ given as an oracle. Then there exists a strongly
 polynomial-time, (in $n$), oracle algorithm to check if $\sum_{1\le i_1
 <\ldots<i_m \le n} \frac{\partial ^m q}{\partial
  x_{i_1}\ldots\partial x_{i_m}} (\0) > 0$.
 \end{corol}

 \proof It follows from Theorem \ref{Tver-new} that $\sum_{1\le i_1
 <\ldots<i_m \le n} \frac{\partial ^m q}{\partial
  x_{i_1}\ldots\partial x_{i_m}} (\0) > 0$ if and only if
  $\frac{\partial^n p}{\partial x_1...\partial x_n} (\0) > 0$ ,
  where $p(\x) = q(\x)(\frac{\sum_{1 \leq i \leq n}
  x_{i}}{n})^{n-m}$ is positive hyperbolic polynomial
  of degree $n$ . Notice the polynomial $q$ is easy to evaluate given
  an oracle evaluating the polynomial $p$ .
  Let $q(\x) = \sum_{ \sum_{1 \leq i \leq n} r_i = n}
  a_{r_{1},...,r_{n}} x_{1}^{r_{1}}...x_{n}^{r_{n}}$ ,
  the support $Supp(q) = \{(r_{1},...,r_{n}) : a_{r_{1},...,r_{n}} \neq 0$ ,
  the Newton polytope is the convex hull
  of the support $CO(Supp(q))$ . It was proved in \cite{Gur-Hyp}
  that an integer vector $(r_{1},...,r_{n}) \in Supp(q)$
  if and only if $(r_{1},...,r_{n}) \in CO(Supp(q))$ .
  Corollary 4.3 in \cite{Gur1}
  provides a strongly polynomial (in $n$)
  oracle algorithm for the membership problem $"X \in CO(Supp(q)) ?"$
  for positive hyperbolic polynomial
  $q(\x)$ of degree $n$ .

  \qed

 \begin{corol}\label{2cor}
  Let $q( \x)$ be a positive hyperbolic polynomial of degree
  $m < n$ given as an oracle. Then
  \begin{enumerate}
  \item
  There exists a strongly
  polynomial-time, (in $n$), oracle algorithm to check if

  $D_{m}(q):=\sum_{1\le i_1
  <\ldots<i_m \le n} \frac{\partial ^m q}{\partial
   x_{i_1}\ldots\partial x_{i_m}} (\0)  > 0$.
  \item
  There exists a deterministic
  polynomial-time oracle algorithm which computes a number $F_{m}(q)$
  satisfying the inequality
  $1 \leq \frac{F_{m}(q)}{D_{m}(q)} \leq 2 \gamma(n,m) \leq 2
  \frac{n^{n}}{n!}$.
  \end{enumerate}
  \end{corol}

 \begin{corol}\label{perm-m}
 Let $A$ be $n \times n$ matrix with nonnegative entries. Then there exists a deterministic
  polynomial-time algorithm which computes a number $P_{m}(A)$
  satisfying the inequality
 $ 1 \leq \frac{P_{m}(A)}{\per_{m} A} \leq 2 \gamma(n,m) \leq 2
 \frac{n^{n}}{n!}$.
 \end{corol}
 It s very possible that there exists a deterministic
 polynomial-time algorithm (in $n$) which approximates $\per_{m} A $ within multiplicative factor
 $e^m$ .

 \subsection{A conjecture}

 Let the assumptions of Corollary \ref{2cor} holds.
 To approximate $D_m(q)$ we used the identity
 \begin{equation}\label{identr}
 \frac{\partial^n}{\partial x_1...\partial x_n} q(\0)r(\0) =
 D_{q}(m),
 \end{equation}
 for $r(\x)=((n-m)!)^{-1}(\sum_{i=1}^n x_i)^{n-m}$.
 It is natural to ask we can improve our estimates if we choose
 a different positive hyperbolic $r(\x)$ such that (\ref{identr})
 holds.

  \begin{defn}
  Denote the set of positive hyperbolic polynomials of degree $m$ and in $n$ variables as $PHYP(n,m)$;
  and by $UPHYP(n,m)\subset PHYP(n,m)$ the subset of polynomials satisfying
  $ \frac{\partial ^{n-m}}{\partial
 x_{i_1}\ldots\partial x_{i_{n-m}}} r(\0)=1$ for all $1 \leq i_1 < \ldots< i_{n-m} \leq
 n$.  Define
  $$
  \gamma(n,m) = \inf_{r \in UPHYP(n,n-m)} \sup_{q \in PHYP(n,m)} \frac{ \Cp qr}
  {\frac{\partial^n}{\partial x_1...\partial x_n}
  q(0,...0)r(0,...,0)}
  .$$
  \end{defn}

  Assume that $r\in  UPHYP(n,n-m)$. Since $r$ is positive hyperbolic,
  all its monomials have nonnegative coefficients.  Hence $r(\x)\ge
  S_{n,n-m}(\x)$ for any $\x\ge 0$.  In particular $\Cp r\ge \Cp
  S_{n,n-m}={n \choose m}$.  By choosing in $\inf\sup$ definition $\gamma(n,m)$
 $q = (\frac{x_1 +...+x_{n}}{n})^{m}$ we deduce straightforward that
 $\gamma(n,m)\ge \frac{ n^{m} }{m!}$.

 \begin{con}\label{new-van}
 $$\gamma(n,m) = \sup_{q \in PHYP(n,m)} \frac{ \Cp q S_{n , n-m}}
 {\frac{\partial^n}{\partial x_1...\partial x_n} q(\0) S_{n ,n-m}(\0)} = \frac{ n^{m}
 }{m!}.$$
 \end{con}
 Note that the hyperbolic van der Waerden inequality \cite{Gur1} implies that
 $\gamma(n,m) \leq \frac{n^{n}}{n!}$.
 It follows from Theorem \ref{frtver}
 that for each $q\in PHYP(n,m)$ such that $\Cp q=q(\1)$ we have the
 inequality
 $$\frac{ \Cp q S_{n , n-m}}
 {\frac{\partial^n}{\partial x_1...\partial x_n} q(\0) S_{n ,n-m}(\0)} \le \frac{ n^{m}
 }{m!}.$$
 \begin{rem}
 We presented in this section one "natural" generating polynomial for
 the hafnian, and described all symmetric boolean matrices
 such that this polynomial is positive hyperbolic. It is an interesting
 open problem whether there exists a generating positive
 hyperbolic polynomial for the hafnians of boolean matrices which can can
 evaluated in polynomial time. If the answer to this problem is negative
 it can explain why approximating the hafnian (number of perfect matchings
 in general graphs) is "harder" than the same problem for the permanent.
 It is easy to prove that computing the hafnian of integer symmetric
 $2n \times 2n$ matrices with nonnegative entries and the signature
 $(+,-,...,-)$
 if $\#P$-complete. The results in this paper allow to use
 Sinkhorn's scaling to approximate the hafnian within multiplicative
 factor $e^{2n}$ in this "hyperbolic" case.
 \end{rem}

\end{document}